\documentstyle[11pt]{amsart}

\newtheorem{thm}{Theorem}[section]
\newtheorem{cor}[thm]{Corollary}
\newtheorem{prop}[thm]{Proposition}
\newtheorem{lem}[thm]{Lemma}

\theoremstyle{definition}
\newtheorem{dfn}[thm]{Definition}

\theoremstyle{remark}
\newtheorem{rem}[thm]{Remark}
\newtheorem{exa}[thm]{Example}
\newtheorem{question}[thm]{Question}

\numberwithin{equation}{section}

\newcommand{\iso}{\stackrel{\simeq}{\rightarrow}}
\newcommand{\inj}{\hookrightarrow}
\newcommand{\surj}{\mbox{$\rightarrow \! \! \! \! \! \rightarrow$}}

\newcommand{\ar}{\rightarrow}
\newcommand{\mx}{{\frak m}}
\newcommand{\opn}{\operatorname}
\renewcommand{\mathrm}[1]{\text{\rom{\rm #1}}}
\renewcommand{\mathsf}[1]{\text{\rom{\sf #1}}}
\renewcommand{\mathbf}[1]{\text{\rom{\bf #1}}}
\newcommand{\bdot}{{\textstyle \cdot}}
\newcommand{\dodot}{*}
\newcommand{\ul}{\underline}
\newcommand{\blnk}[1]{\mbox{\hspace{#1}}}
\newcommand{\rmitem}[1]{\item[\text{\rom{(#1)}}]}
\newcommand{\exar}[1]{ @> #1 >> }
\newcommand{\sqbr}[1]{ [\, #1 \, ] }

\title[Beilinson Completion Algebras]{Traces and Differential
Operators over Beilinson Completion Algebras}
\author[Amnon Yekutieli]{Amnon Yekutieli*}
\address{Department of Theoretical Mathematics,
The Weizmann Institute of Science,
Rehovot 76100, ISRAEL}
\date{8 August 1994}
\email{amnon@@wisdom.weizmann.ac.il}
\thanks{* Supported by an Alon Fellowship, and incumbent
of the Maurice M.\ Boukstein Career Development Chair.}

\newcommand{\lrar}[1]{\begin{picture}(50,10)(-25,-5)
\put(-25,0){\vector(1,0){50}}
\put(0,5){\makebox(0,0)[b]{\mbox{$#1$}}}
\end{picture}}

\newcommand{\ldar}[1]{\begin{picture}(10,50)(-5,-25)
\put(0,25){\vector(0,-1){50}}
\put(5,0){\mbox{$#1$}}
\end{picture}}

\newcommand{\luar}[1]{\begin{picture}(10,50)(-5,-25)
\put(0,-25){\vector(0,1){50}}
\put(5,0){\mbox{$#1$}}
\end{picture}}

\newcommand{\ldrar}[1]{\begin{picture}(50,50)(-25,-25)
\put(-25,25){\vector(1,-1){50}}
\put(5,0){\mbox{$#1$}}
\end{picture}}

\begin{document}
\setcounter{section}{-1}
\maketitle


\section{Introduction}

In the short paper \cite{Be} A.\ Beilinson introduced a generalized
version of adeles, with values in any quasi-coherent sheaf on a noetherian
scheme $X$. In particular, taking the structure sheaf $\cal{O}_{X}$
one gets the cosimplicial ring of adeles
$\Bbb{A}^{\bdot}(X, \cal{O}_{X})$. In each degree $n$,
$\Bbb{A}^{n}(X, \cal{O}_{X})$ is a subring (a ``restricted product'')
of the product of local factors
$\prod_{\xi} \cal{O}_{X, \xi}$. Here $\xi = (x_{0}, \ldots, x_{n})$
runs over all chains of length $n$ of points in $X$.
The Beilinson completion $\cal{O}_{X, \xi}$ is gotten by a process of
inverse and direct limits. For $n=0$, $\cal{O}_{X,(x_{0})}$
is simply the $\frak{m}$-adic completion of the local ring at $x_{0}$.
For applications to duality theory one is primarily interested in the
completion $\cal{O}_{X, \xi}$ along a {\em saturated chain} $\xi$.
As shown in \cite{Ye1}, the
semi-local ring $\cal{O}_{X, \xi}$ carries a natural topology, and
its residue fields carry rank $n$ valuations.

In the present paper we isolate the completion $\cal{O}_{X, \xi}$ from its
geometric environment,
and study it as a separate algebraic-topological object, which we call
a {\em Beilinson completion algebra} (BCA).
The methods used here belong to commutative algebra,
analysis and differential geometry. Our main results
have to do with {\em dual modules} of BCAs, their functorial behavior
and their interaction with differential operators. These results, in turn,
have some noteworthy applications to algebraic geometry (see \S 0.3).

One may view our paper partly as a continuation of the work of Lipman,
Kunz and others on explicit formulations of duality theory
(cf.\ \cite{Li1}, \cite{Li2}, \cite{Ku}, \cite{Hu1}, \cite{Hu2},
\cite{HK1}, \cite{HK2},
\cite{HS}, \cite{LS}, \cite{Hg}). Their work deals with {\em linear} aspects
of duality theory - construction of dualizing modules, trace maps etc.
To that we have little new to add in the present paper. The novelty of our
work is in establishing the {\em nonlinear} properties of duality theory.
We show how duality interacts with {\em differential phenomena}, such as
$\cal{D}$-modules and De Rham complexes.
Such results seem to have been beyond the reach of the methods
of commutative algebra used henceforth in this area.

In the remainder of the introduction we outline the content of the paper.

\subsection{Beilinson completion algebras}
Let $k$ be a fixed perfect base field.
A local BCA $A$ is a quotient of a ring
$F((\ul{s})) [\sqbr{ \ul{t} }] =
F(( s_{1}, \ldots, s_{m} )) [\sqbr{ t_{1}, \ldots, t_{n} }]$,
where $F$ is a finitely generated field extension of $k$, and
$F(( s_{1}, \ldots, s_{m} )) = F((s_{m})) \cdots  ((s_{1}))$
is an iterated field of Laurent series. $A$ is a complete noetherian local
ring, and a
semi-topological (ST) $k$-algebra. On the residue field $A / \frak{m}$
there is a structure of $m$-dimensional topological local field (TLF).
(These terms are explained briefly in \S 1-2.) The surjection
$F((\ul{s})) [\sqbr{ \ul{t} }] \surj A$ is not part of the structure of $A$.
A general BCA is a finite product of local ones.

We are interested in two kinds of homomorphisms between BCAs. The first is
called a {\em morphism of BCAs}, and the second is called an {\em
intensification homomorphism}. Rather than defining these notions here
(this is done in \S 2-3), we demonstrate them by examples.
Let
$A := k(s)[\sqbr{ t }]$ and
$B := k(s)((t))$. These local BCAs arise geometrically: take
$X := \mathbf{A}^{2}_{k} = \opn{Spec} k \sqbr{ s,t }$ and
$x = (0), y = (t), z = (s,t) \in X$.
Then
$A \cong \cal{O}_{X, (y)}$ and
$B \cong \cal{O}_{X, (x, y)}$, the Beilinson completions of $\cal{O}_{X}$
along the chains $(y), (x,y)$ respectively. The inclusion
$A \ar B$ is a morphism, which in ``cosimplicial'' notation is
$\partial^{+} : \cal{O}_{X, (y)} \ar \cal{O}_{X, (x, y)}$.
Now let
$\hat{A} := k((s))[\sqbr{ t }] \cong \cal{O}_{X, (y, z)}$. Then
$A \ar \hat{A}$ is an intensification homomorphism, which we also write
as
$\partial^{-} : \cal{O}_{X, (y)} \ar \cal{O}_{X, (y, z)}$.

Whenever $A \ar B$ is a morphism and $A \ar \hat{A}$ is an intensification,
there is a BCA
$\hat{B} = B \otimes_{A}^{(\wedge)} \hat{A}$, a morphism
$\hat{A} \ar \hat{B}$ and an intensification $B \ar \hat{B}$. This
situation is called {\em intensification base change}. In our example,
$\hat{B} = k((s))((t)) \cong \cal{O}_{X, (x, y, z)}$.

BCAs and morphisms of BCAs constitute a category which is denoted by
$\mathsf{BCA}(k)$.

\subsection{The Results}
There are three main results in the paper. Their precise statement is in
the body of the paper, and what follows is only a sketch.

A {\em finite type ST module} $M$ over a BCA $A$ is a quotient of
$A^{n}$ for some $n$, with the quotient topology (so if $A / \frak{m}$ is
discrete, $M$ has the $\frak{m}$-adic topology.) The {\em fine topology}
on an $A$-module $M$ is characterized by the property that each finitely
generated submodule $M' \subset M$, with the subspace topology,
is of finite type.
(More on ST modules in \S 1.)
Given a TLF $K$ (i.e.\ a BCA which is a field), we denote by
$\omega(K)$ the top degree component of the separated algebra of
differentials
$\Omega_{K/k}^{\bdot, \mathrm{sep}}$.

\bigskip \noindent
{\bf Theorem \ref{thm6.1}:}\ (Dual modules)\ Let $A$ be a local BCA and $M$
a finite type ST $A$-module. Then there is a {\em dual module}
$\opn{Dual}_{A} M$, enjoying the following properties. To any morphism
$\sigma: K \ar A$ in $\mathsf{BCA}(k)$
with $K$ a field, there is a bijection
\[ \Psi_{\sigma}^{M} : \opn{Dual}_{A} M \iso
\opn{Hom}_{K; \sigma}^{\mathrm{cont}}(M, \omega(K)).  \]
If $\sigma = \tau \circ f$ for some morphisms $f : K \ar L$ and
$\tau : L \ar A$, then
\[ \Psi_{\sigma}^{M}(\phi) = \opn{Res}_{L/K} \circ \Psi_{\tau}^{M}(\phi), \]
where
$\opn{Res}_{L/K} : \omega(L) \ar \omega(K)$ is the residue on TLFs, see
\cite{Ye1} \S 2.4. If
$\sigma, \sigma' : K \ar A$ are two pseudo coefficients fields (i.e.\
morphisms such that
$[A / \frak{m} : K] < \infty$) which are congruent modulo $\frak{m}$, then
the isomorphism
\[ \Psi_{\sigma, \sigma'}^{M} = \Psi_{\sigma'}^{M} \circ
(\Psi_{\sigma}^{M})^{-1} :
\opn{Hom}_{K; \sigma}^{\mathrm{cont}}(M, \omega(K)) \iso
\opn{Hom}_{K; \sigma'}^{\mathrm{cont}}(M, \omega(K)) \]
has an explicit formula in terms of ``Taylor expansions'' and differential
operators.

\bigskip
In particular for $M=A$ we set
$\cal{K}(A) := \opn{Dual}_{A} A$, with the fine topology.
$\cal{K}(A)$ is an injective hull of the residue field $A / \frak{m}$.
Note that for a field
$K$, $\cal{K}(K) = \omega(K)$.
If $M$ is any ST $A$-module we define
\[ \opn{Dual}_{A} M := \opn{Hom}_{A}^{\mathrm{cont}}(M, \cal{K}(A)) \]
with the $\opn{Hom}$ topology. (When $M$ is of finite type this is consistent
with Thm.\ \ref{thm6.1}.)
We show that given an intensification homomorphism
$v : A \ar \hat{A}$ there is a continuous homomorphism of ST $A$-modules
\[ q_{\hat{A} / A}^{M} = q_{v}^{M} : \opn{Dual}_{A} M \ar
\opn{Dual}_{\hat{A}} (\hat{A} \otimes_{A} M). \]

\bigskip \noindent
{\bf Theorem \ref{thm7.2}:}\ (Traces)\
Let $A \ar B$ be a morphism in $\mathsf{BCA}(k)$. Then there exists a
continuous $A$-linear trace map
$\opn{Tr}_{B/A} : \cal{K}(B) \ar \cal{K}(A)$. This trace is functorial:
$\opn{Tr}_{C/A} = \opn{Tr}_{B/A} \circ \opn{Tr}_{C/B}$. It induces a
bijection
\[ \cal{K}(B) \iso \opn{Hom}_{A}^{\mathrm{cont}}(B, \cal{K}(A)). \]
The trace
commutes with intensification base change: given an intensification
$A \ar \hat{A}$, and letting
$\hat{B} := B \otimes_{A}^{(\wedge)} \hat{A}$, we have
\[ q_{\hat{A} / A} \circ \opn{Tr}_{B / A} =
\opn{Tr}_{\hat{B} / \hat{A}} \circ q_{\hat{B} / B}. \]
If $\sigma : K \ar A$ is a morphism with $K$ a field, then
\[ \opn{Tr}_{A / K}(\phi) = \Psi_{\sigma}^{A}(\phi)(1) \in \omega(K) \]
for $\phi \in \cal{K}(A)$.

\bigskip \noindent
{\bf Theorem \ref{thm8.1}:}\ (Duals of continuous differential operators)\
Suppose $M, N$ are ST $A$-modules with the fine topologies and
$D : M \ar N$ is a continuous DO. Then there is a continuous DO
\[ \opn{Dual}_{A}(D) : \opn{Dual}_{A} N \ar \opn{Dual}_{A} M. \]
This operation is transitive in $D$ and compatible with intensification
base change $A \ar \hat{A}$. $\opn{Dual}_{A}(D)$ is unique, has an explicit
description using the isomorphisms
$\Psi_{\sigma}^{M}, \Psi_{\sigma}^{N}$, and is the adjoint of
$D$ w.r.t.\ suitably defined residue pairings.

\subsection{Applications}
The primary application of our results, and the original motivation of
the paper, is the explicit construction of residue complexes on $k$-schemes.
This is carried out in \cite{Ye2}. The construction is extremely simple,
and we shall sketch it here. Suppose $X$ is a $k$-scheme of finite type
and $(x, y)$ is a saturated chain of points in it (i.e.\ $y$ is an immediate
specialization of $x$). There are natural homomorphisms
$\partial^{-} : \cal{O}_{X, (x)} \ar \cal{O}_{X, (x, y)}$ and
$\partial^{+} : \cal{O}_{X, (y)} \ar \cal{O}_{X, (x, y)}$, the first being
an intensification and the second a morphism (cf.\ example in \S 0.1 above).
According to Theorems \ref{thm6.1} and \ref{thm7.2} we get an
$\cal{O}_{X}$-linear homomorphism
\[ \delta_{(x, y)} : \cal{K}(\cal{O}_{X, (x)})
\exar{q_{\partial^{-}}} \cal{K}(\cal{O}_{X, (x, y)})
\exar{\opn{Tr}_{\partial^{+}}} \cal{K}(\cal{O}_{X, (y)}). \]
Considering $\cal{K}(\cal{O}_{X, (x)})$ as a skyscraper sheaf sitting on
$\{ x \}^{-}$, we define
\begin{eqnarray*}
\cal{K}^{\bdot}_{X} & := & \bigoplus_{x \in X} \cal{K}(\cal{O}_{X, (x)}) \\
\delta_{X} & := & \sum_{(x, y)} \delta_{(x, y)}.
\end{eqnarray*}
Then $(\cal{K}^{\bdot}_{X}, \delta_{X})$ is the residue complex on $X$
(cf.\ \cite{RD}, \cite{EZ}, \cite{Ye1} and \cite{SY}).

A special feature of this particular construction of $\cal{K}^{\bdot}_{X}$
is that given a DO
$D : \cal{M} \ar \cal{N}$ between $\cal{O}_{X}$-modules, there is a dual DO
\begin{equation} \label{eqn0.1}
\opn{Dual}_{X}(D) : \cal{H}om_{\cal{O}_{X}}(\cal{N}, \cal{K}^{\bdot}_{X})
\ar \cal{H}om_{\cal{O}_{X}}(\cal{M}, \cal{K}^{\bdot}_{X})
\end{equation}
which is a homomorphism of complexes.
This implies that $\cal{K}^{\bdot}_{X}$ is a complex of right
$\cal{D}_{X}$-modules. Conversely, $\cal{D}_{X}$ can be recovered from DOs
acting on $\cal{K}^{\bdot}_{X}$.
Another consequence of (\ref{eqn0.1}) is that
$\cal{F}^{\bdot \bdot}_{X} :=
\cal{H}om_{\cal{O}_{X}}(\Omega^{\bdot}_{X/k}, \cal{K}^{\bdot}_{X})$
has a natural structure of double complex.
Using $\cal{F}^{\bdot \bdot}_{X}$ we are able to analyze the niveau
spectral sequence converging to
$\mathrm{H}^{\mathrm{DR}}_{\bdot}(X)$, the algebraic De Rham homology of $X$.

\subsection{Plan of the paper} \blnk{1mm} \\
{\bf Section 1}: a quick review of semi-topological rings and modules, as
well as new facts on ST $\opn{Hom}$ modules.

\medskip \noindent
{\bf Section 2}: definition of BCAs and morphisms, including examples.

\medskip \noindent
{\bf Section 3}: definition of intensification homomorphisms, base change.

\medskip \noindent
{\bf Section 4}: general facts on continuous differential operators over ST
algebras; the Lie derivative.

\medskip \noindent
{\bf Section 5}: the structure of the ring of continuous DOs $\cal{D}(K)$
over a TLF $K$; $\omega(K)$ is a right $\cal{D}(K)$-module, and
the action is by adjunction in a suitable sense.

\medskip \noindent
{\bf Section 6}: existence of dual modules is proved.

\medskip \noindent
{\bf Section 7}: contravariance of dual modules w.r.t.\ morphisms is
proved (traces).

\medskip \noindent
{\bf Section 8}: the interaction between dual modules and DOs is examined,
leading to Thm. \ref{thm8.1} and a few corollaries.

\medskip \noindent
{\bf Acknowledgements.}\
I wish to thank J.\ Lipman for his continued interest in this work.
Also thanks to V.\ Lunts for very helpful remarks on differential
operators.


\section{Some Results on Semi-Topological Rings}

Let us recall some definitions and results from \cite{Ye1} \S 1.
A semi-topological (ST) ring is a ring $A$, with a linear topology on
its underlying additive group, such that
for all $a \in A$, left and right multiplication by $a$ are continuous maps
$\lambda_{a}, \rho_{a} : A \ar A$. A
ST left $A$-module is an $A$-module $M$, whose underlying additive group is
linearly topologized, and such that for all $a \in A$ and $x \in M$,
the multiplication maps they define
$\lambda_{a} : M \ar M$ and $\rho_{x} : A \ar M$ are continuous. ST left
$A$-modules and continuous $A$-linear homomorphisms form a category, denoted
$\mathsf{STMod}(A)$. Similarly one defines ST right modules and bimodules.

Assume for simplicity that the ST ring $A$ is commutative. In
$\mathsf{STMod}(A)$ there are direct and inverse limits, and a tensor
product. Given a ST $A$-module $M$, the associated separated module
$M^{\opn{sep}} = M / \{ 0 \}^{-}$ is also a ST $A$-module.
The category $\mathsf{STMod}(A)$ is additive, but not abelian. An exact
sequence in it is, by definition, a sequence
$M' \exar{ \phi } M \exar{ \psi } M''$ which is exact in the untopologized
sense (i.e.\ in $\mathsf{Mod}(A)$), and such that both $\phi$ and $\psi$
are strict.

On any $A$-module
$M$ there is a finest topology making it into a ST module; it is called the
fine $A$-module topology. If $M$ has the fine topology, then for any ST
$A$-module $N$, one has
$\opn{Hom}_{A}^{\opn{cont}}(M,N) = \opn{Hom}_{A}(M,N)$,
and this in fact characterizes the fine topology.
Trivially, if $M$ has the fine topology, then so does $M^{\mathrm{sep}}$.
A free ST $A$-module is a free $A$-module with the fine topology. So $F$ is
free iff
$F \cong \bigoplus A$ with the $\bigoplus$ topology. A ST module $M$ has
the fine topology iff it admits a strict surjection $F \surj M$ with $F$
free.

\begin{dfn}
Let $M,N$ be ST $A$-modules. The (weak) $\opn{Hom}$ topology on the abelian
group
$\opn{Hom}_{A}^{\opn{cont}}(M,N)$ is the coarsest linear topology such that
for every $x \in M$, the map
$\rho_{x} : \opn{Hom}_{A}^{\opn{cont}}(M,N) \ar N$, $\phi \mapsto \phi(x)$,
is continuous.
\end{dfn}

Unless otherwise specified, this is the topology we consider on
$\opn{Hom}_{A}^{\opn{cont}}(M,N)$. If $M$ has the fine topology, we shall
often drop the superscript ``$\opn{cont}$''.

\begin{rem}
A basis of neighborhoods of $0$ for the $\opn{Hom}$ topology is the
collection of open subgroups
$\{ V(F,U) \}$, where $F$ runs over the finite subsets of $M$, $U$ runs over
the open subgroups of $N$, and
$V(F,U) = \{ \phi\ |\ \phi(F) \subset U \}$. Such a topology is sometimes
called the weak topology (cf.\ \cite{Ko}). Usually, to obtain a duality one
needs a finer topology - the strong topology of \cite{Ko}, or the
compact-open topology of \cite{Mc}. In the present
paper duality is defined by indirect means, and for our purposes the
weak topology suffices (cf.\ Remark \ref{rem8.1}).
\end{rem}

The next lemma summarizes the properties of the $\opn{Hom}$ topology. Its
easy proof is left to the reader.

\begin{lem} \label{lem1.1}
Let $A$ be a commutative ST ring.
\begin{enumerate}
\item Let $\phi : M' \ar M$ and $\psi : N \ar N'$ be homomorphisms in
$\mathsf{STMod}(A)$. Then the induced homomorphism
$\opn{Hom}_{A}^{\opn{cont}}(M,N) \ar \opn{Hom}_{A}^{\opn{cont}}(M',N')$
is continuous.

\item Let $M,N$ be ST $A$-modules. Then
$\opn{Hom}_{A}^{\opn{cont}}(M,N)$ is a ST $A$-module. \linebreak
$\opn{End}_{A}^{\opn{cont}}(M) = \opn{Hom}_{A}^{\opn{cont}}(M,M)$ is a ST
$A$-algebra, and $M$ is a ST left
$\opn{End}_{A}^{\opn{cont}}(M)$-module.
The natural bijection
$M \iso \opn{Hom}_{A}^{\opn{cont}}(A,M)$, $x \mapsto \rho_{x}$, is
an isomorphism of ST $A$-modules.

\item Suppose in \rom{(1)} $\phi$ is surjective and $\psi$ is a strict
monomorphism. Then
$\opn{Hom}_{A}^{\opn{cont}}(M,N) \ar \opn{Hom}_{A}^{\opn{cont}}(M',N')$
is a strict monomorphism.

\item Let $(M_{\alpha})_{\alpha \in I}$ be a direct system in
$\mathsf{STMod}(A)$, with $I$ a directed set. Then for any ST $A$-module $N$
the natural map
\[ \lim_{\leftarrow \alpha} \opn{Hom}_{A}^{\opn{cont}}(M_{\alpha},N) \ar
\opn{Hom}_{A}^{\opn{cont}}(\lim_{\alpha \ar} M_{\alpha}, N) \]
is an isomorphism of ST $A$-modules.
\end{enumerate}
\end{lem}

{}From parts (1) and (2) of the lemma it follows that
$\opn{Hom}_{A}^{\opn{cont}}$ is an
additive bifunctor
$\mathsf{STMod}(A)^{\circ} \times \mathsf{STMod}(A) \ar \mathsf{STMod}(A)$.

Tensor products are defined in $\mathsf{STMod}(A)$. The usual tensor
product $M \otimes_{A} N$ is given the finest linear topology s.t.\
the maps
$\rho_{y}: M \ar M \otimes_{A} N$, $x' \mapsto x' \otimes y$ and
$\lambda_{x}: N \ar M \otimes_{A} N$, $y' \mapsto x \otimes y'$ are all
continuous (see \cite{Ye1} Def.\ 1.2.11).

\begin{lem} \label{lem1.2}
\rom{(Adjunction)}\ Let $A,B$ be ST rings (not necessarily commutative), let
$L$ be a ST left $A$-module, $N$ a ST left $B$-module, and $M$ a ST
$B$-$A$-bimodule. Then
\[ \opn{Hom}_{B}^{\opn{cont}}(M \otimes_{A} L, N) \cong
\opn{Hom}_{A}^{\opn{cont}}(L, \opn{Hom}_{B}^{\opn{cont}}(M, N)) \]
as topological abelian groups.
\end{lem}

\begin{pf}
Immediate from the definitions of the $\opn{Hom}$ and $\otimes$
topologies.
\end{pf}

We say a homomorphism
$\phi: M \ar N$ of topological abelian groups is dense if
$\phi(M) \subset N$ is (everywhere) dense.

\begin{lem} \label{lem1.3}
Suppose $A$ is a ST ring and $M \ar \hat{M}$, $N \ar \hat{N}$ are
continuous dense homomorphisms of ST $A$-modules. Then
$M \otimes_{A} N \ar \hat{M} \otimes_{A} \hat{N}$ is dense.
\end{lem}

\begin{pf}
By transitivity of denseness it suffices to prove that
$M \otimes_{A} N \ar \hat{M} \otimes_{A} N$ is dense. Choose a surjection
from a free module
$A^{(I)} = \bigoplus A$ onto $N$. This induces surjections
$M^{(I)} \surj M \otimes_{A} N$ and
$\hat{M}^{(I)} \surj \hat{M} \otimes_{A} N$. But according to \cite{Ye1}
Prop.\ 1.1.8 (c), $M^{(I)} \ar \hat{M}^{(I)}$ is dense.
\end{pf}

\begin{dfn} \label{def1.1}
Let $A$ be a commutative noetherian ST ring. A ST $A$-module $M$ is called
of {\em finite type}
(resp.\ {\em cofinite type}, resp.\ {\em torsion type}) if it is
finitely generated (resp.\ it is artinian, resp.\
$\opn{Supp} M \subset \opn{Spec} A$ consists solely of maximal ideals),
and if it has the fine topology.
\end{dfn}

Denote the
full subcategories of $\mathsf{STMod}(A)$ consisting of finite type
(resp.\ cofinite type) modules by
$\mathsf{STMod}_{\opn{f}}(A)$ (resp.\ $\mathsf{STMod}_{\opn{cof}}(A)$).

Generalizing the Zariski and Artin-Rees properties for noetherian rings with
adic topologies, we make the following definition. Let us point out that
this definition is stronger then \cite{Ye1} Definition 3.2.10.

\begin{dfn} \label{def1.2}
Let $A$ be a noetherian commutative ST ring. $A$ is said to be a {\em
Zariski ST ring} if
\begin{enumerate}
\rmitem{i} Every ST $A$-module, which is either of finite type or of
torsion type, is separated.

\rmitem{ii} Every (continuous) $A$-linear homomorphism between two ST
$A$-modules, each either of finite type or of torsion type, is strict.
\end{enumerate}
\end{dfn}

\begin{prop} \label{prop1.1}
Let $A$ be a local Zariski ST ring, with maximal ideal $\mx$. Assume that
$A \cong \lim_{\leftarrow i} A / \mx^{i+1}$ as ST rings. Let $M,N$ be
ST $A$-modules.
\begin{enumerate}
\item If $M,N$ are both of finite type then so is
$\opn{Hom}_{A}^{\mathrm{cont}}(M,N)$.

\item If $M$ is of finite type and $N$ is of cofinite type then
$\opn{Hom}_{A}^{\mathrm{cont}}(M,N)$ is of cofinite type.

\item If $M,N$ are both of cofinite type then
$\opn{Hom}_{A}^{\mathrm{cont}}(M,N)$ is of finite type.
\end{enumerate}
\end{prop}

\begin{pf}
(1)\ Let $A^{r} \surj M$ be a surjection. By Lemma \ref{lem1.1}
(2) and (3),
$\opn{Hom}_{A}^{\mathrm{cont}}(M,N)$ \linebreak
$\inj N^{r}$
is a strict monomorphism. Now use the
Zariski property to conclude that $\opn{Hom}_{A}(M,N)$ has the fine
topology.

\medskip \noindent (2)\ Like (1).

\medskip \noindent (3)\
Let $M_{i} := \opn{Hom}_{A}^{\mathrm{cont}}(A / \mx_{i+1} ,M)$, so
$M_{i} \inj M$ is strict, $M_{i}$ has the fine topology, and
$M = \lim_{i \ar} M_{i}$. Similarly define $N_{i}$.
By part (4) of Lemma \ref{lem1.1},
\[ \opn{Hom}_{A}^{\mathrm{cont}}(M,N) =
\lim_{\leftarrow i} \opn{Hom}_{A}^{\mathrm{cont}}(M_{i}, N) =
\lim_{\leftarrow i} \opn{Hom}_{A}^{\mathrm{cont}}(M_{i}, N_{i}). \]
Now $M_{i}$ and
$N_{i}$ are of finite type, so we can we can use part (1) and
\cite{Ye1} Prop.\ 1.2.20.
\end{pf}

\begin{cor} \label{cor1.1}
\rom{(ST Version of Matlis Duality)}\
Let $A$ be as in the proposition. Suppose
$I$ is an injective hull of $A/\mx$, endowed with the fine topology. Then
$\opn{Hom}_{A}^{\mathrm{cont}}(-,I)$ is an equivalence
\[ \mathsf{STMod}_{\opn{f}}(A)^{\circ} \leftrightarrow
\mathsf{STMod}_{\opn{cof}}(A)\ . \]
\end{cor}


\section{Definitions and Basic Properties of BCAs}

In this section $k$ is a fixed perfect field.
If $A$ is a ST $k$-algebra and
$\ul{t} = (t_{1}, \ldots, t_{n})$ is a sequence of indeterminates, we
denote by
$A [\sqbr{\ul{t}}] = A [\sqbr{t_{1}, \ldots, t_{n}}]$
the ring of formal power series, with the topology given by
\[ A [\sqbr{\ul{t}}] = \lim_{\leftarrow i} A \sqbr{\ul{t}} / (\ul{t})^{i},
\]
where for each $i$, $A \sqbr{ \ul{t} } / (\ul{t})^{i}$ has the fine
$A$-module topology.
The ring of Laurent series $A((t))$ is topologized by
\[ A((t)) = \lim_{j \ar} t^{-j} A [\sqbr{t}], \]
and we define recursively
\[ A((\ul{t})) = A((t_{1}, \ldots, t_{n})) :=
A((t_{2}, \ldots, t_{n}))((t_{1})). \]
According to \cite{Ye1} \S 1.3,
$A [\sqbr{\ul{t}}]$ and $A((\ul{t}))$ are ST $k$-algebras.

A {\em topological local field} (TLF) over $k$ is a field $K$, together with
a topology, and valuation rings $\cal{O}_{i}$, $i=1, \ldots, n$,
such that the residue field $\kappa_{i}$ of $\cal{O}_{i}$ is the fraction
field of $\cal{O}_{i+1}$, and
$K = \opn{Frac}(\cal{O}_{1})$. These data are related by the existence of
a {\em parametrization}: an isomorphism
$K \cong F(( t_{1}, \ldots, t_{n} ))$ of ST $k$-algebras, s.t.\
$\cal{O}_{i} \cong F(( t_{i+1}, \ldots, t_{n} ))
[\sqbr{ t_{i} }]$. Here $F$ is a discrete field, and
$\Omega^{1}_{F/k}$ has finite rank. The number $n$ is the dimension of
the local field $K$.
Topological local fields constitute a category $\mathsf{TLF}(k)$.
For more details see \cite{Ye1} \S 2.1.

\begin{dfn} \label{dfn2.1}
A local {\em Beilinson completion algebra} (BCA) over $k$
is a commutative semi-topological local ring
$A$, together with a structure of topological local field on the
residue field $A/\mx$. The following
condition must be satisfied: there exists a surjective homomorphism of
$k$-algebras
\[ F(( \ul{s} )) [\sqbr{ \ul{t} }] =
F(( s_{1}, \ldots, s_{m} )) [\sqbr{ t_{1}, \ldots, t_{n} }] \ar A, \]
which is strict (topologically), and induces and isomorphism
of TLFs $F(( \ul{s} )) \cong A/ \mx$. Such a surjection is called a
parametrization of $A$.

A Beilinson completion algebra is a finite product of local BCAs.
\end{dfn}

\begin{rem}
In greater generality one can define a BCA over any noetherian ring $R$,
to be any finite algebra over the $R$-algebra
$\Bbb{A}(\Xi, \cal{O}_{X}) = \prod_{\xi \in \Xi} \cal{O}_{X, \xi}$,
where $\Xi$ is a finite set of saturated chains in some finite type
$R$-scheme $X$, and $\Bbb{A}(- , -)$ is Beilinson's scheme theoretical
group of adeles. See \cite{Be}, \cite{Hr}, \cite{Ye1} and \cite{HY} for
the definition of adeles, and cf.\ Examples \ref{exa2.0} and \ref{exa2.1}
below.
\end{rem}

Observe that a Beilinson completion algebra $A$ is necessarily an
$\frak{r}$-adically complete, noetherian, semi-local ring, where
$\frak{r}$ is the Jacobson radical of $A$.
If $A$ is artinian, then in the terminology of \cite{Ye1}, it is a
cluster of TLFs (a CTLF).

For any $\mx \in \opn{Max} A$
set $\opn{res.dim}_{\mx} A := \opn{dim} A/\mx$, the local field dimension.
We say that $A$ is equidimensional of dimension $n$ if
$\opn{res.dim}_{\mx} A =n$ for all $\mx$. In this case  we set
$\opn{res.dim} A := n$, and
\begin{eqnarray*}
\cal{O}_{i}(A) & := & \prod_{\mx \in \opn{Max} A} \cal{O}_{i}(A/\mx) \\
\kappa_{i}(A) & := & \prod_{\mx \in \opn{Max} A} \kappa_{i}(A/\mx) \\
\end{eqnarray*}
for $1 \leq i \leq n$. Also we set
$\cal{O}_{0}(A) := A$ and
$\kappa_{0}(A) := A / \frak{r} = \prod A / \mx$.

The motivating example is:

\begin{exa} \label{exa2.0}
Let $X$ be a scheme of finite type over $k$, and let
$\xi = (x_{0}, \ldots, x_{m})$ be a saturated chain in $X$. Then the
Beilinson completion $\cal{O}_{X, \xi}$ of the structure sheaf along $\xi$
is defined; see \cite{Ye1} \S 3.1.
We claim that $\cal{O}_{X, \xi}$ is an equidimensional BCA, of dimension
$m$. To see why, first choose a coefficient field
$\sigma: k(x_{0}) \ar \hat{\cal{O}}_{X, x_{0}} = \cal{O}_{X, (x_{0})}$.
According to \cite{Ye1} Lemma 3.3.9, $\sigma$ extends to a lifting
$\sigma_{\xi} : k(\xi) = k(x_{0})_{\xi} \ar \cal{O}_{X, \xi}$. Sending
$t_{1}, \ldots, t_{n}$ to generators of the maximal ideal
$\frak{m}_{x_{0}}$,  we get a strict surjection
$k(\xi) [\sqbr{ t_{1}, \ldots, t_{n} }] \surj \cal{O}_{X, \xi}$. Finally,
according to \cite{Ye1} Prop.\ 3.3.6, $k(\xi)$ is a finite product of
TLFs, all of dimension $m$.
\end{exa}

\begin{exa} \label{exa2.1}
Consider a BCA
$A = F(( s_{1}, \ldots, s_{m} )) [\sqbr{ t_{1}, \ldots, t_{n} }]$. We
claim it is of the form $\cal{O}_{X, \xi}$. Choose an
integral $k$-scheme of finite type $Y$ such that $F = k(Y)$. Set
$X:= \mathbf{A}_{Y}^{n+m} = \mathbf{A}_{k}^{n+m} \times_{k} Y$, and let
$\xi=(x_{0}, \ldots, x_{m})$ be the saturated chain
$x_{i} := (t_{1}, \ldots, t_{n}, s_{1}, \ldots, s_{i})$, where
we write
$\mathbf{A}_{k}^{n+m} = \opn{Spec} k[\ul{s},\ul{t}]$. Then
$F((\ul{s}))[[\ul{t}]] \cong \cal{O}_{X,\xi}$ (cf.\ \cite{Ye1} Thm.\ 3.3.2
(c); it can be assumed that $Y$ is normal).
\end{exa}

Let $A$ be a local BCA of $\opn{res.dim}$ $n$. For every $1 \leq i \leq n$
there is a subring
$\cal{O}_{1, \ldots, i}(A) \subset A$ defined by
\[ \cal{O}_{1, \ldots, i}(A) \cong
A \times_{A / \frak{m}}
\cal{O}_{1}(A / \frak{m}) \times_{\kappa_{1}(A / \frak{m})} \cdots
\times_{\kappa_{i-1}(A / \frak{m})} \cal{O}_{i}(A / \frak{m}).  \]
It is the largest subring of $A$ which projects onto $\kappa_{i}(A)$,
and it is actually the valuation ring of a rank $i$ valuation (hence
local). In \cite{Ye1} the notation
$\cal{O}(A)$ was used for $\cal{O}_{1, \ldots, n}(A)$.

\begin{dfn} \label{dfn2.2} (Morphisms)\
Let $A$ and $B$ be Beilinson completion algebras. A morphism $f: A \ar B$
is a continuous $k$-algebra homomorphism, satisfying the following local
condition.
Given a maximal ideal $\frak{n} \subset B$, let $\mx \subset A$ be the
unique maximal ideal such that $f^{-1}(\frak{n}) \subset \mx$. Set
$i := \opn{res.dim} B_{\frak{n}} - \opn{res.dim} A_{\frak{m}}$,
which is assumed to be non-negative. Then
$f(A_{\frak{m}}) \subset \cal{O}_{1, \ldots, i}(B_{\frak{n}})$,
the induced homomorphism
$A_{\frak{m}} \ar \kappa_{i}(B/\frak{n})$ sends $\frak{m}$ to $0$,
and
$A/\mx \ar \kappa_{i}(B/\frak{n})$ is a finite morphism of local
fields.
\end{dfn}

The composition of two morphisms is again a morphism, so we get a category,
which is denoted by $\mathsf{BCA}(k)$.
The number $i$ in the definition is called the relative residual dimension of
$f$ at $\frak{n}$, denoted $\opn{res.dim}_{\frak{n}} f$. If $f$ is
equidimensional we shall ommit the subscript $\frak{n}$. We call $f$ finite
if $B$ is a finitely generated $A$-module.
Observe that the full subcategory of $\mathsf{BCA}(k)$ consisting of fields
coincides with the full subcategory of $\mathsf{TLF}(k)$ consisting of
TLFs whose last residue field is finitely generated over $k$. (In
characteristic $0$ this is all of $\mathsf{TLF}(k)$.)

Here are some typical examples of morphisms of BCAs.

\begin{exa}
Let $A:= k[\sqbr{s}]$, $B:= k((s))[\sqbr{t}]$, and let $f: A \ar B$ be the
inclusion. Then $\mx=(s)$, $\frak{n} = (t)$, $\opn{res.dim}_{\mx} A =0$,
$\opn{res.dim}_{\frak{n}} B =1$ and
$\opn{res.dim}_{\frak{n}} f = 1$.
\end{exa}

\begin{exa} \label{exa2.2}
Let $X$ be a finite type $k$-scheme, $\xi=(x,\ldots,y)$ a saturated chain in
$X$, $A:= \cal{O}_{X,(y)}$, $B:= \cal{O}_{X,\xi}$, and
$\partial^{+} : \cal{O}_{X,(y)} \ar \cal{O}_{X,\xi}$ the coface map. Now
$\opn{res.dim} A = 0$, and
$\opn{res.dim} B = \opn{res.dim} \partial^{+}$ equals the length of $\xi$.
\end{exa}

\begin{exa}
Let $X,Y$ be finite type $k$-schemes, $f:X \ar Y$ a
$k$-morphism, $y \in Y$ any point and $x$ a closed point in the fibre
$X_{y} := f^{-1}(y)$. Since $k(y) \ar k(x)$ is finite,
$f^{*}:  \cal{O}_{Y,(y)} \ar \cal{O}_{X,(x)}$ is a morphism of BCAs, with
$\opn{res.dim} f^{*} = 0$.
\end{exa}

\begin{dfn}
Let $A$ be a local BCA over $k$, with maximal ideal $\mx$. A
{\em coefficient field} (resp.\ {\em quasi coefficient field}, resp.\
{\em pseudo coefficient field}) for $A$ is a morphism
$\sigma: K \ar A$ in $\mathsf{BCA}(k)$, with $K$ a field, and such that
the induced homomorphism
$K \ar A/\mx$ is bijective (resp.\ finite separable, resp.\ finite).
\end{dfn}

By definition, every local BCA has a coefficient field.

\begin{lem} \label{lem2.2}
Let $A$ be a local BCA over $k$, with maximal ideal $\mx$. Then:
\begin{enumerate}
\rmitem{a} Suppose $A$ is artinian and $K \ar A$ is a pseudo coefficient
field. Then $A$ has the fine $K$-module topology.

\rmitem{b} Letting $A_{i} := A / \mx^{i+1}$,
the map $A \ar \lim_{\leftarrow i} A_{i}$ is an isomorphism of ST
$k$-algebras.

\rmitem{c} Let $K \ar A$ be a pseudo coefficient field, and let $M$ be a
torsion type ST $A$-module (see Def.\ \ref{def1.1}). Then $M$ is a free
ST $K$-module.

\rmitem{d} Suppose $\sigma: K \ar A$ is a morphism of BCAs, with $K$ a
field. Then there exists a finite morphism
$f: L [\sqbr{ \ul{t} }] \ar A$ extending $\sigma$, i.e.\
$\sigma: K \ar L \ar L [\sqbr{ \ul{t} }] \exar{f} A$.
\end{enumerate}
\end{lem}

\begin{pf}
(a)\ By \cite{Ye1} Prop.\ 2.2.2.

\medskip \noindent(b)\ This is true for $F((\ul{s}))[\sqbr{\ul{t}}]$
(by definition!) and hence, by \cite{Ye1} Prop.\ 1.2.20,
for every quotient $A$.

\medskip \noindent (c)\ Set
$M_{i} := \opn{Hom}_{A}(A_{i}, M)$, with the fine $A$-module topology.
According to \cite{Ye1} Cor.\ 1.2.6, $M \cong \lim_{i \ar} M_{i}$. Now
$M_{i}$ is a ST $A_{i}$-module with the
fine topology. Since $A_{i}$ has the fine $K$-module topology, so does
$M_{i}$. Passing to the limit, $M$ has the fine $K$-module topology, so it
is a free ST $K$-module.

\medskip \noindent (d)\
According to \cite{Ye1} Cor.\ 2.1.19 we can find a finite morphism
$K(( \ul{s} )) = L \ar A / \frak{m}$. As in the proof of ibid.\ Prop.\
2.2.2., this extends to a morphism
$L \ar \lim_{\leftarrow i} A_{i} = A$, which we then extend to
$f: L [\sqbr{ \ul{t} }] \ar A$ by sending the $t_{i}$ to generators of
the maximal ideal ideal $\frak{m}$.
\end{pf}

\begin{prop} \label{prop2.1}
Let $A$ be a BCA over $k$. Then:
\begin{enumerate}
\rmitem{a} If $f : A \ar B$ is a finite morphism in $\mathsf{BCA}(k)$,
then $B$ has the fine $A$-module topology.

\rmitem{b} Conversely, if $B$ is a finite $A$-algebra, then $B$ admits
a unique structure of BCA s.t.\ $A \ar B$ is a morphism of BCAs.

\rmitem{c} $A$ is a Zariski ST ring. Moreover, every finite type or torsion
type ST $A$-module is complete.
\end{enumerate}
\end{prop}

\begin{pf}
(a)\
Let $\frak{r} \subset A$ and $\frak{s} \subset B$ be the Jacobson radicals.
According to \cite{Ye1} Prop.\ 2.2.2 (b),
$B_{i} := B/ \frak{s}^{i+1}$ has the fine
$A_{i} := A / \frak{r}^{i+1}$-module topology, for
each $i \geq 0$. So $B_{i}$ also has the fine $A$-module
topology. Now use Lemma \ref{lem2.2} (b) and \cite{Ye1} Prop.\ 1.2.20.

\medskip \noindent (b)\ According to \cite{Ye1} Prop.\ 2.2.2 (c), this is
true for $A_{i} \ar B_{i}$. Now use
$B \cong \lim_{\leftarrow i} B_{i}$.

\medskip \noindent (c)\
It suffices to consider $A = F((\ul{s}))[\sqbr{\ul{t}}]$.
By \cite{Ye1} Thm.\ 3.3.8, $A$ is a Zariski ST ring
in the sense of ibid.\ Def.\ 3.2.10. This means that every finite type ST
$A$-module is separated, and every homomorphism between two such modules
is strict.

Now consider two torsion type ST $A$-modules, $M$ and $N$. We may assume $A$
is local. Choose a pseudo coefficient field $K \ar A$. Then $M,N$ are free ST
$K$-modules, and in particular they are separated and complete (cf.\
\cite{Ye1} Prop.\ 1.5). To prove that any homomorphism $\phi : M \ar N$
is strict, we may assume it is injective. Then any $K$-linear splitting
$M \surj N$ is continuous, showing that $\phi$ is strict.

Finally, given a homomorphism $\phi : M \ar N$, with $M,N$ either of
finite type or of torsion type, then the module
$\bar{M} := \phi(M)$, endowed with the fine topology, is a ST module of both
types. Therefore $M \surj \bar{M}$ and $\bar{M} \inj N$ are both strict.
\end{pf}


\section{Intensification Base Change}

The operation of base change to be discussed in this subsection
is a generalization of the one in \cite{Ye1} \S 2.2.
The important notion is that of an intensification homomorphism
$u : A \ar \hat{A}$ between two BCAs (Def.\ \ref{dfn3.2}). Differentially
$u$ is ``\'{e}tale``: the differential invariants of $\hat{A}$ descend
to $A$. From the point of view of valuations, $\hat{A}$ is like a completion
of $A$.
Again $k$ is a fixed perfect field.

\begin{dfn} \label{dfn3.1}
Let $A,\hat{A} \in \mathsf{BCA}(k)$ have Jacobson radicals
$\frak{r},\hat{\frak{r}}$
respectively, and let $u: A \ar \hat{A}$ be a continuous $k$-algebra
homomorphism, with $u(\frak{r}) \subset \hat{\frak{r}}$.
\begin{enumerate}
\rmitem{a} $u$ is called {\em radically unramified} if
$\hat{\frak{r}} = \hat{A} \cdot u(\frak{r})$.

\rmitem{b} $u$ is called {\em finitely ramified} if
$\hat{A} / \hat{A} \cdot u(\frak{r})$ is artinian, and if for every
$\hat{\mx} \in \opn{Max} \hat{A}$ lying over some
$\mx \in \opn{Max} A$, letting $n := \opn{res.dim} \hat{A} / \hat{\frak{m}}$,
the image of $(A/\mx)^{\times}$ in the rank $n$ valuation group of
$\hat{A} / \hat{\mx}$
has finite index.
\end{enumerate}
\end{dfn}

\begin{prop} \label{prop3.1}
\rom{(Finitely Ramified Base Change)}\
Let $K, \hat{K}, A \in \mathsf{BCA}(k)$, with $A$ a local ring and
$K, \hat{K}$ fields. Suppose $f : K \ar A$ is a morphism in
$\mathsf{BCA}(k)$ and $u: K \ar \hat{K}$ is a finitely ramified homomorphism.
Then there exists a BCA $\hat{A}$, a morphism
$\hat{f} : \hat{K} \ar \hat{A}$ in $\mathsf{BCA}(k)$, and a finitely
ramified homomorphism $v: A \ar \hat{A}$, satisfying:
\begin{enumerate}
\rmitem{i} $v \circ f = \hat{f} \circ u$, and moreover the homomorphism
$A \otimes_{K} \hat{K} \ar \hat{A}$ is dense.

\rmitem{ii} $\opn{res.dim} \hat{f} = \opn{res.dim} f$.

\rmitem{iii} Suppose $\hat{g} : \hat{K} \ar \hat{C}$ is a morphism in
$\mathsf{BCA}(k)$, with $\hat{C}$ local, and let
$n:= \opn{res.dim} \hat{g} - \opn{res.dim} \hat{f}$. Suppose also
$w: A \ar \hat{C}$ is a continuous homomorphism s.t.\
$w \circ f = \hat{g} \circ u$,
$w(A) \subset \cal{O}_{1, \ldots, n}(\hat{C})$, and
$A \ar \kappa_{n}(\hat{C})$ is finitely ramified.
Then there exists a unique morphism
$\hat{h} : \hat{A} \ar \hat{C}$ (of $\opn{res.dim}$ $n$)
in $\mathsf{BCA}(k)$, such that
$\hat{g}  = \hat{h} \circ \hat{f}$ and
$w =  \hat{h} \circ v$.
\end{enumerate}
\end{prop}

\begin{pf}
Choose a finite morphism
$K(( \ul{s} ))[\sqbr{ \ul{t} }] \ar A$ (cf.\ Lemma \ref{lem2.2}), and set
\[ \hat{A} := A \otimes_{K(( \ul{s} ))[\sqbr{ \ul{t} }]}
\hat{K}(( \ul{s} ))[\sqbr{ \ul{t} }]. \]
$\hat{A}$ is a BCA by Prop.\ \ref{prop2.1}, and $\hat{f}, v$ are the
obvious maps.

Let us prove that $A \otimes_{K} \hat{K} \ar \hat{A}$ is dense.
Denoting by
$K \sqbr{ \ul{s}, \ul{s}^{-1} }$ the ring of Laurent polynomials,
we have
$\hat{K} \otimes_{K} A \cong
\hat{K} \sqbr{ \ul{s}, \ul{s}^{-1}, \ul{t} }
\otimes_{K \sqbr{ \ul{s}, \ul{s}^{-1}, \ul{t} }} A$.
By \cite{Ye1} Lemma 1.3.9 the homomorphism
$\hat{K} \sqbr{ \ul{s}, \ul{s}^{-1} } \ar \hat{K}(( \ul{s} ))$ is dense,
and a similar argument shows that so is
$\hat{K} \sqbr{ \ul{s}, \ul{s}^{-1}, \ul{t} } \ar
\hat{K}(( \ul{s} ))[\sqbr{ \ul{t} }]$. Now use Lemma \ref{lem1.3}.

Finally, given $\hat{C}$, the arguments in the proof of \cite{Ye1}
Thm.\ 2.2.4 imply there is a morphism
$\hat{K}(( \ul{s} ))[\sqbr{ \ul{t} }] \ar \hat{C}$, and tensoring with
$A$ we get
$\hat{h}: \hat{A} \ar \hat{C}$. Uniqueness follows from the denseness of
$A \otimes_{K} \hat{K} \ar \hat{A}$ .
\end{pf}

The algebra $\hat{A}$ in the proposition is unique (up to a unique
isomorphism). We shall denote it by
\begin{equation} \label{eqn3.2}
\hat{A} = A \otimes^{(\wedge)}_{K} \hat{K}\ .
\end{equation}
In contrast with the usual tensor product, this is not a symmetric
expression - we shall always put the algebra which is the range of the
finitely ramified homomorphism to the right.

In \cite{Ye1} \S 1.5 the notion of a topologically \'{e}tale homomorphism
relative to $k$ was defined. A homomorphism
$v: A \ar \hat{A}$ in $\mathsf{STComAlg}(k)$, the category of commutative ST
$k$-algebras, is called topologically \'{e}tale
relative to $k$ if for any separated ST $\hat{A}$-module $\hat{M}$, any
continuous $k$-linear derivation $\partial : A \ar \hat{M}$ has a
unique extension to a continuous derivation
$\hat{\partial} : \hat{A} \ar \hat{M}$. Often we shall suppress the phrase
``relative to $k$''; this should not cause any confusion as we have no
notion of absolute topologically \'{e}tale homomorphism.

\begin{lem} \label{lem3.1} \mbox{ }
\begin{enumerate}
\rmitem{a} The homomorphism
$v: A \ar \hat{A} = A \otimes^{(\wedge)}_{K} \hat{K}$ is flat.
\rmitem{b} If $u : K \ar \hat{K}$ is topologically \'{e}tale relative to $k$,
then $v: A \ar \hat{A}$ is topologically \'{e}tale  and radically
unramified.
\end{enumerate}
\end{lem}

\begin{pf}
(a)\ We have
$\hat{A} \cong A \otimes_{K((\ul{s}))[\sqbr{\ul{t}}]}
\hat{K}((\ul{s}))[\sqbr{\ul{t}}]$.
According to \cite{CA} Ch.\ III \S 5.4 Prop.\ 4, the homomorphism
$K((\ul{s}))[\sqbr{\ul{t}}] \ar \hat{K}((\ul{s}))[\sqbr{\ul{t}}]$
is flat; hence so is $A \ar \hat{A}$.

\medskip \noindent (b)\
As in the proof of \cite{Ye1} Thm.\ 2.4.23,
$K((\ul{s}))[\sqbr{\ul{t}}] \ar \hat{K}((\ul{s}))[\sqbr{\ul{t}}]$
is topologically \'{e}tale. By \cite{Ye1} Prop.\ 1.5.9 (b), so is
$A \ar \hat{A}$. The ring
$\hat{A} / \hat{A} \cdot v(\mx) \cong A / \mx \otimes_{K((\ul{s}))}
\hat{K}((\ul{s}))$
is  reduced, since
$K((\ul{s})) \ar \hat{K}((\ul{s}))$ is separable (cf.\ proof of \cite{Ye1}
Thm.\ 2.4.23). This shows that $\hat{A} \cdot v(\mx)$ is the Jacobson
radical of $\hat{A}$.
\end{pf}

Let $A,\hat{A}$ be two local BCAs, with maximal ideals $\mx,\hat{\mx}$
respectively. Suppose $v: A \ar \hat{A}$ is a finitely ramified, radically
unramified homomorphism. Let $\sigma : K \ar A$ be a pseudo coefficient field,
and assume there is some subfield
$\hat{K} \subset \hat{A} / \hat{\mx}$ such that
$K \ar \hat{K}$ is topologically \'{e}tale relative to $k$, and
$A/ \mx \otimes_{K} \hat{K} \ar \hat{A} / \hat{\mx}$ is bijective.
Then $\hat{K} \ar \hat{A} / \hat{\mx}$ is finite, and
$K \ar \hat{K}$ is finitely ramified. Also, this $\hat{K}$ is unique.
Since some
lifting $\hat{K} \ar \hat{A}$ exists, there is a unique pseudo coefficient
field
\begin{equation} \label{eqn3.1}
\hat{\sigma} : \hat{K} \ar \hat{A}
\end{equation}
extending $\sigma$ (cf.\ \cite{Ye1} formula (4.1.11)).

\begin{exa} \label{exa3.1}
If $v: A \ar \hat{A}$ is topologically \'{e}tale and
$K \ar A/\mx$ is purely inseparable, then such a subfield $\hat{K}$ exists.
Indeed, we have
$\hat{A} / \hat{\mx} = \hat{A} \otimes_{A} A/\mx$, so
$A / \mx \ar \hat{A} / \hat{\mx}$ is also topologically
\'{e}tale. If $\sigma$ is a coefficient field, the statement is trivial.
Otherwise, see \cite{Ye1} formula (4.1.10).
\end{exa}

We make
$\opn{gr}_{\mx} A = \bigoplus_{i \geq 0} \mx^{i} / \mx^{i+1}$
into a graded ST ring by putting on
$\mx^{i} \subset A$ the subspace topology, and putting on
$\mx^{i} / \mx^{i+1}$ the quotient topology. Similarly, for a ST $A$-module
$M$, $\opn{gr}_{\mx} M$ is a graded ST $\opn{gr}_{\mx} A$-module.

\begin{prop} \label{prop3.2}
In the situation above, suppose that $v : A \ar \hat{A}$ is flat. Then:
\begin{enumerate}
\rmitem{a} For any finite type ST $A$-module $M$ which has finite length,
the canonical homomorphism
\[ \hat{K} \otimes_{K} M \ar  \hat{A} \otimes_{A} M \]
is an isomorphism of ST $\hat{K}$-modules.
\rmitem{b} The canonical morphism
$A \otimes^{(\wedge)}_{K} \hat{K} \ar \hat{A}$ in $\mathsf{BCA}(k)$
is an isomorphism.
\rmitem{c} For any finite type ST $A$-module $M$,
the canonical homomorphism
\[ \hat{K} \otimes_{K} \opn{gr}_{\mx} M \ar  \opn{gr}_{\mx} (\hat{A}
\otimes_{A} M) \]
is an isomorphism of graded ST $\hat{K}$-modules.
\end{enumerate}
\end{prop}

\begin{pf}
(a)\ The proof is by induction on the length of $M$. For $M$ of length $1$,
we have by assumption
\[ \hat{K} \otimes_{K} M \cong \hat{K} \otimes_{K} A/\mx \cong
\hat{A} / \hat{\mx} \cong \hat{A} \otimes_{A} M\ . \]
Otherwise, we can find an exact sequence (of untopologized $A$-modules)
\[ M^{\bdot} = ( 0 \ar M' \ar M \ar M'' \ar 0 ) \]
which gives rise, by flatness, to a homomorphism of exact sequences
$\hat{K} \otimes_{K} M^{\bdot} \ar \hat{A} \otimes_{A} M^{\bdot}$. By
induction and the Five Lemma, we conclude that
$\hat{K} \otimes_{K} M \cong \hat{A} \otimes_{A} M$. Since both modules have
the fine $\hat{K}$-module topologies, this is a homeomorphism.

\medskip \noindent (b)\
We have
$A \otimes_{K}^{(\wedge)} \hat{K} \cong \lim_{\leftarrow i}
A / \mx^{i+1} \otimes_{K} \hat{K}$, and by Lemma \ref{lem2.2} (b),
$\hat{A} \cong \lim_{\leftarrow i} \hat{A} /  \hat{\mx}^{i+1}$.
Now use part (a) above, together with the isomorphism
$\hat{A} \otimes_{A} (A / \mx^{i+1}) \cong \hat{A} /  \hat{\mx}^{i+1}$.

\medskip \noindent (c)\
By flatness and the fact that $A$ and $\hat{A}$ are Zariski ST rings, it
follows that
$\hat{A} \otimes_{A} \mx^{i} M \cong \mx^{i} (\hat{A} \otimes_{A} M)
\subset \hat{A} \otimes_{A} M$ as ST $\hat{A}$-modules. Therefore
$\hat{A} \otimes_{A} (\opn{gr}_{\mx} M)_{i} \cong
\opn{gr}_{\mx} (\hat{A} \otimes_{A} M)_{i}$
as ST $\hat{A} / \hat{\mx}$-modules. Now use part (a).
\end{pf}

\begin{dfn} \label{dfn3.2} \rom{(Intensification)}\
Let $u: A \ar \hat{A}$ be a continuous $k$-algebra homomorphism between
two BCAs. If $u$ is flat, finitely ramified, radically unramified and
topologically \'{e}tale relative to $k$, then $u$ is called an
{\em intensification} homomorphism.
\end{dfn}

\begin{exa} \label{exa3.2}
Let $X$ be a finite type $k$-scheme, $\xi=(x,\ldots,y)$ a saturated chain in
$X$, $A:= \cal{O}_{X,(x)}$, $B:= \cal{O}_{X,\xi}$, and
$\partial^{-} : \cal{O}_{X,(x)} \ar \cal{O}_{X,\xi}$ the coface map. Then
$\partial^{-}$ is an intensification homomorphism (cf.\ Example
\ref{exa2.2}).
\end{exa}

\begin{thm} \label{thm3.1} \rom{(Intensification Base Change)}\
Let $A, \hat{A}, B$ be local BCAs, let $f : A \ar B$ be a morphism
in $\mathsf{BCA}(k)$, and let $u : A \ar \hat{A}$ be an intensification
homomorphism. Then there is a BCA
$\hat{B} = B \otimes_{A}^{(\wedge)} \hat{A}$, a morphism
$\hat{f} :  \hat{A} \ar \hat{B}$ and an intensification homomorphism
$v : B \ar \hat{B}$, satisfying conditions \rom{(i) - (iii)} of Prop.\
\ref{prop3.1} \rom{(}but replacing the letters $K,A$ with $A,B$\rom{)}.
\end{thm}

\medskip \noindent
\begin{equation} \label{eqn3.3}
\setlength{\unitlength}{0.25mm}
\begin{array}{ccc}
B & \lrar{v} 	& \hat{B} \\
\luar{f} &	& \luar{\hat{f}} \\
A & \lrar{u} 	& \hat{A}
\end{array}
\end{equation}

\begin{pf}
Choose a coefficient field $\sigma : K \ar A$, and let
$\hat{\sigma} : \hat{K} = K \otimes_{A} \hat{A} \ar \hat{A}$ be its unique
extension. So $\hat{A} \cong A \otimes_{K}^{(\wedge)} \hat{K}$. Set
$\hat{B} := B \otimes_{K}^{(\wedge)} \hat{K}$. We can find a surjective
morphism
$K[\sqbr{\ul{t}}] \surj A$, and it gives
$\hat{A} \cong A \otimes_{K[\sqbr{\ul{t}}]} \hat{K}[\sqbr{\ul{t}}] $.
The homomorphism
$\hat{K}\sqbr{\ul{t}} \ar B$ extends uniquely to a morphism
$\hat{K}[\sqbr{\ul{t}}] \ar B$ : define it inductively into
$\cal{O}_{i}(B)$, $i= \opn{res.dim} f, \ldots, 2, 1$. Hence
$\hat{f} : \hat{A} \ar \hat{B}$ is also defined. The uniqueness of $\hat{B}$
is clear from its construction.
\end{pf}

\begin{exa}
Let
$A := k(s)[\sqbr{ t }]$,
$\hat{A} := k((s))[\sqbr{ t }]$ and
$B := k(s)((t))$, so the inclusion $A \ar \hat{A}$
(resp.\ $A \ar B$) is an intensification (resp.\ a morphism). We then have
\[ k(s)((t)) \otimes_{k(s)[\sqbr{ t }]}^{(\wedge)} k((s))[\sqbr{ t }]
\cong k((s))((t)). \]
\end{exa}

\begin{prop} \label{prop3.3} \rom{(Associativity)}\
Say
$C \leftarrow B \ar \Hat{B} \leftarrow \Hat{A} \ar \Hat{\Hat{A}}$
are BCAs and homomorphisms, where the ``$\leftarrow$'' are morphisms,
and the ``$\ar$'' are intensifications. Then there is a canonical
isomorphism of BCAs
\[ (C \otimes_{B}^{(\wedge)} \Hat{B}) \otimes_{\Hat{A}}^{(\wedge)}
\Hat{\Hat{A}} \cong
C \otimes_{B}^{(\wedge)}
(\Hat{B} \otimes_{\Hat{A}}^{(\wedge)} \Hat{\Hat{A}}). \]
\end{prop}

\begin{pf}
Set
$\Hat{\Hat{B}} := \Hat{B} \otimes_{\Hat{A}}^{(\wedge)} \Hat{\Hat{A}}$
and
$\Hat{C} := C \otimes_{B}^{(\wedge)} \Hat{B}$.
By construction (cf.\ Prop.\ \ref{prop3.1}) we get an intensification
homomorphism
$\Hat{C} \ar C \otimes_{B}^{(\wedge)} \Hat{\Hat{B}}$, and together with
the morphism
$\Hat{\Hat{A}} \ar \Hat{\Hat{B}} \ar
C \otimes_{B}^{(\wedge)} \Hat{\Hat{B}}$
we deduce, using Cor. \ref{thm3.1}, the existence of a morphism
$h: \Hat{C} \otimes_{\Hat{A}}^{(\wedge)} \Hat{\Hat{A}} \ar
C \otimes_{B}^{(\wedge)} \Hat{\Hat{B}}$.
The same corollary says there is a morphism
$\Hat{\Hat{B}} \ar \Hat{C} \otimes_{\Hat{A}}^{(\wedge)} \Hat{\Hat{A}}$,
and together with the intensification
$C \ar \Hat{C} \otimes_{\Hat{A}}^{(\wedge)} \Hat{\Hat{A}}$
we get a morphism
$h': C \otimes_{B}^{(\wedge)} \Hat{\Hat{B}} \ar
\Hat{C} \otimes_{\Hat{A}}^{(\wedge)} \Hat{\Hat{A}}$. Since the maps
from
$C \otimes_{B} \Hat{B} \otimes_{\Hat{A}} \Hat{\Hat{A}}$
to both these BCAs are dense, $h$ and $h'$ must be each other's inverse.
\end{pf}


\section{Continuous Differential Operators}

We begin with some general results on continuous differential
operators (DOs) over ST algebras. Let $k$ be a discrete commutative ring,
let $A$
be a commutative, separated, ST $k$-algebra, and let $M$ be a separated ST
$A$-module. For $n \geq 0$, the separated module of principal parts
$\cal{P}_{A/k}^{n,\opn{sep}}$ is the ST $A$-$A$-bimodule
$(A \otimes_{k} A / I^{n+1})^{\opn{sep}}$, where
$I := \opn{ker}(A \otimes_{k} A \ar A)$.
Set
$\cal{P}_{A/k}^{n,\opn{sep}}(M) :=
(\cal{P}_{A/k}^{n,\opn{sep}} \otimes_{A} M)^{\opn{sep}}$, which is an
$A$-module by
$a \cdot (( 1 \otimes 1) \otimes x) = (a \otimes 1) \otimes x$.
The universal continuous DO of order $n$ is
$\opn{d}_{M}^{n} : M \ar \cal{P}_{A/k}^{n,\opn{sep}}(M)$,
$\opn{d}_{M}^{n}(x) = (1 \otimes 1) \otimes x$
(see \cite{EGA} Ch.\ IV \S 16.8 and \cite{Ye1} \S 1.5).
For any separated ST $A$-module $N$, $\opn{d}_{M}^{n}$ induces a bijection
\[ \opn{Hom}_{A}^{\opn{cont}}(\cal{P}_{A/k}^{n,\opn{sep}}(M),N) \cong
\opn{Diff}_{A/k}^{n,\opn{cont}}(M,N)\ . \]
There are inclusions
\[ \opn{Diff}_{A/k}^{n,\opn{cont}}(M,N) \subset
\opn{Diff}_{A/k}^{\opn{cont}}(M,N) \subset
\opn{Hom}_{k}^{\opn{cont}}(M,N)\ . \]
$\opn{Diff}_{A/k}^{\opn{cont}}(M,N)$ is a filtered $A$-$A$-bimodule, where for
$D \in \opn{Diff}_{A/k}^{\opn{cont}}(M,N)$ and $a,b \in A$ we have
$a D b = a \circ D \circ b : M \ar N$.
Denote the order of the DO $D$ by $\opn{ord}_{A}(D)$.

\begin{dfn} \label{4.1}
Given a separated ST $A$-module $M$, let
$\cal{D}(A;M) := \opn{Diff}_{A/k}^{\opn{cont}}(M,$ \linebreak
$ M)$,
which is a filtered $k$-algebra.
For $M=A$ we shall write simply
$\cal{D}(A) := \cal{D}(A;A)$.
\end{dfn}

Denote the left action of $\cal{D}(A;M)$ on $M$ by $D \dodot x$,
for $D \in \cal{D}(A;M)$ and $x \in M$.

\begin{rem}
$\cal{D}(A)$ can be made into a ST $k$-algebra by giving it the subspace
topology w.r.t. the embedding
$\cal{D}(A) \subset \opn{End}_{k}^{\opn{cont}}(A)$.
However we shall not make use of this  topology.
\end{rem}

\begin{lem} \label{lem4.3}
Assume that for some $n \geq 0$,
$\cal{P}_{A/k}^{n,\opn{sep}}$ is a finite type ST left $A$-module. If
$M$ is a finite type ST $A$-module, then so is
$\cal{P}_{A/k}^{n,\opn{sep}}(M)$.
\end{lem}

\begin{pf}
First note that
$\cal{P}_{A/k}^{n,\opn{sep}}$ is a commutative ST ring, admitting two
continuous $k$-algebra homomorphisms
$A \ar \cal{P}_{A/k}^{n,\opn{sep}}$. By \cite{Ye1} Cor.\ 4.5,
$\cal{P}_{A/k}^{n,\opn{sep}} \otimes_{A} M$ is a finite type ST
$\cal{P}_{A/k}^{n,\opn{sep}}$-module. The left $A$-module structure on
$\cal{P}_{A/k}^{n,\opn{sep}}$ comes from the algebra homomorphism
$a \mapsto a \otimes 1$. From \cite{Ye1} Prop.\ 2.9 and our
assumption it follows that
$\cal{P}_{A/k}^{n,\opn{sep}} \otimes_{A} M$ has the fine $A$-module
topology. But then the same is true for
$\cal{P}_{A/k}^{n,\opn{sep}}(M) =
(\cal{P}_{A/k}^{n,\opn{sep}} \otimes_{A} M)^{\opn{sep}}$.
\end{pf}

Define
\[ \cal{T}(A) := \opn{Der}^{\opn{cont}}_{k}(A,A) = \opn{Hom}^{\opn{cont}}_{A}
(\Omega^{1,\opn{sep}}_{A/k},A)\ . \]
Corresponding to the decomposition
$\cal{P}_{A/k}^{1,\opn{sep}} = A \oplus \Omega^{1,\opn{sep}}_{A/k}$
we have
$A \oplus  \cal{T}(A) = \cal{D}^{1}(A) \subset \cal{D}(A)$,
and just like in the discrete case, $\cal{T}(A)$ is a Lie algebra over $k$.

\begin{lem} \label{lem4.4}
Suppose $\hat{A}$ is another commutative, separated, ST $k$-algebra, and
$u: A \ar \hat{A}$ is a topologically \'{e}tale homomorphism relative to $k$.
Then there is an induced homomorphism of filtered $k$-algebras
$\cal{D}(A) \ar \cal{D}(\hat{A})$, sending an operator
$D : A \ar A$ to its unique extension
$\hat{D} : \hat{A} \ar \hat{A}$. More generally, if $M$ is a ST
$A$-module, there is a homomorphism
$\cal{D}(A; M^{\opn{sep}}) \ar
\cal{D}(\hat{A}; (\hat{A} \otimes_{A} M)^{\opn{sep}})$.
\end{lem}

\begin{pf}
The existence and uniqueness of this ring homomorphism are immediate
consequences of \cite{Ye1} Thm.\ 1.5.11 (iv).
\end{pf}

The ring homomorphism
$\cal{D}(A) \ar \cal{D}(\hat{A})$ restricts to a Lie algebra homomorphism
$\cal{T}(A) \ar \cal{T}(\hat{A})$.

\begin{prop} \label{prop4.1}
Let $A,\hat{A}$ be separated ST $k$-algebras, and
let $u : A \ar \hat{A}$ be a flat, topologically \'{e}tale homomorphism
relative to $k$. Assume that for every $n \geq 0$,
$\cal{P}_{A/k}^{n,\opn{sep}}$ is a finitely presented, finite type ST left
$A$-module.
Then the homomorphism
$\cal{D}(A) \ar \cal{D}(\hat{A})$ induces an isomorphism of
filtered $\hat{A}$-$\cal{D}(A)$-bimodules
\[ \hat{A} \otimes_{A} \cal{D}(A) \iso \cal{D}(\hat{A})\ . \]
\end{prop}

\begin{pf}
Since
$\otimes$ commutes with $\lim_{\ar}$, it suffices
to prove that for all $n \geq 0$,
$\hat{A} \otimes_{A} \cal{D}^{n}(A) \ar \cal{D}^{n}(\hat{A})$
is bijective. The assumptions imply that
\[ \begin{array}{rcl}
\hat{A} \otimes_{A} \cal{D}^{n}(A) & = &
\hat{A} \otimes_{A}
\opn{Hom}^{\opn{cont}}_{A}(\cal{P}_{A/k}^{n,\opn{sep}},A) \\
& = & \opn{Hom}_{\hat{A}}^{\opn{cont}}
(\hat{A} \otimes_{A} \cal{P}_{A/k}^{n,\opn{sep}},\hat{A}) \\
& \cong & \opn{Hom}_{\hat{A}}^{\opn{cont}}
(\cal{P}_{\hat{A}/k}^{n,\opn{sep}},\hat{A}) \\
& = & \cal{D}^{n}(\hat{A}) .
\end{array} \]
\end{pf}

Now consider the separated algebra of differentials
$\Omega^{\bdot,\opn{sep}}_{A/k}$, which is a graded ST $k$-algebra
(see \cite{Ye1} Def.\ 1.5.3).  Then
\[ \cal{T}(\Omega^{\bdot,\opn{sep}}_{A/k}) :=
\opn{Der}^{\opn{cont}}_{k}(\Omega^{\bdot,\opn{sep}}_{A/k},
\Omega^{\bdot,\opn{sep}}_{A/k}) \]
is a graded Lie algebra. For instance, the exterior
derivative $\opn{d}$ is an element of degree $1$ in
$\cal{T}(\Omega^{\bdot,\opn{sep}}_{A/k})$.

We shall need  a version of the Lie derivative for semi-topological
algebras (see \cite{Wa} \S 2.24 for the differentiable manifold version).

\begin{prop} \label{prop4.2} \rom{(Lie derivative)}\
Let $A$ be a separated ST $k$-algebra and let $\partial$ be a
continuous $k$-derivation of $A$. Then there exists a unique
continuous, degree $0$, $k$-linear derivation $\opn{L}_{\partial}$ of
$\Omega^{\bdot,\opn{sep}}_{A/k}$, which extends $\partial$ and commutes with
$\mathrm{d}$. The map $\partial \mapsto {\rm L}_{\partial}$ is a
homomorphism of $k$-Lie algebras
$\cal{T}(A) \ar \cal{T}(\Omega^{\bdot,\opn{sep}}_{A/k})$, and is functorial
with respect to topologically \'{e}tale homomorphisms $A \ar \hat{A}$ in
$\mathsf{STComAlg}(k)$.
\end{prop}

\begin{pf}
Let
$\partial \in \cal{T}(A) = \opn{Der}^{\opn{cont}}_{k}(A,A)$ be given.
Since $A$ is separated we get a continuous $A$-linear map
$\Omega^{1,\opn{sep}}_{A/k} \ar A$, which extends by universality to a
continuous degree $-1$ derivation
$\iota_{\partial}: \Omega^{\bdot,\opn{sep}}_{A/k} \ar
\Omega^{\bdot,\opn{sep}}_{A/k}$,
the interior derivative.

Define $\opn{L}_{\partial} := \iota_{\partial} \circ \mathrm{d} +
\mathrm{d} \circ \iota_{\partial}$ (i.e. the graded commutator of
$\iota_{\partial}$ and $\mathrm{d}$). The properties of
$\opn{L}_{\partial}$
are easily deduced from its definition and the fact that
$\mathrm{d}^{2}=0$. To show uniqueness it suffices to consider
$\opn{L}_{\partial}(\alpha)$ for $\alpha=a$ or $\alpha=\mathrm{d}a$,
$a \in A$.
But $\opn{L}_{\partial}(a)=\partial(a)$ and
$\opn{L}_{\partial}(\mathrm{d}a)= \mathrm{d} (\opn{L}_{\partial}(a))=
\mathrm{d} \circ \partial(a)$.

Now let $\partial_{1},\partial_{2} \in \cal{T}(A)$.
Then $[\opn{L}_{\partial_{1}}, \opn{L}_{\partial_{2}}]$ is a continuous
derivation of $\Omega^{\bdot,\opn{sep}}_{A/k}$ commuting with $\mathrm{d}$,
and for all $a \in A$,
\[ [\opn{L}_{\partial_{1}}, \opn{L}_{\partial_{2}}](a)
=[\partial_{1},\partial_{2}](a)
= \opn{L}_{[\partial_{1},\partial_{2}]}(a)\ , \]
so
$[\opn{L}_{\partial_{1}}, \opn{L}_{\partial_{2}}] =
\opn{L}_{[\partial_{1},\partial_{2}]}$.
The functoriality of $\opn{L}$ follows from the same functoriality of
$\iota$ (and $\mathrm{d}$).
\end{pf}

\begin{lem} \label{lem4.1}
Suppose
$\Omega^{n+1,\opn{sep}}_{A/k}=0$ for some $n$. Then for any
$a \in A$, $\alpha \in \Omega^{n,\opn{sep}}_{A/k}$ and
$\partial \in \cal{T}(A)$, one has
$\opn{L}_{a \partial}(\alpha) = \opn{L}_{\partial}(a \alpha)$.
\end{lem}

\begin{pf}
First note that
\begin{equation} \label{eqn4.6}
\partial(a) \alpha - \mathrm{d}(a) \iota_{\partial}(\alpha) =
\iota_{\partial}(\mathrm{d}(a) \alpha) =  0 .
\end{equation}
Since
$\mathrm{d} \alpha = 0$,
$\iota_{a \partial(\alpha)} = a \iota_{\partial}(\alpha)$ and
$ = \partial(a) \alpha +
a \opn{L}_{\partial}(\alpha)$,
it follows that
$\opn{L}_{a \partial}(\alpha) = \opn{L}_{\partial}(a \alpha)$.
\end{pf}

Now assume $k$ is a perfect field.

\begin{prop} \label{prop4.3}
Let $A \in \mathsf{BCA}(k)$, and let $M$ be a finite type ST $A$-module.
Then for any $n \geq 0$,
$\cal{P}_{A/k}^{n,\opn{sep}}(M)$ is a finite type ST left $A$-module.
In particular, so is
$\cal{P}_{A/k}^{1,\opn{sep}} = A \oplus \Omega^{1,\opn{sep}}_{A/k}$, so $A$
is differentially of finite type over $k$ \rom{(}In the sense of
 \cite{Ye1} Def.\ \rom{1.5.16)}.
\end{prop}

\begin{pf}
We may assume $A$ is a local ring. Choose a parametrization of $A$, i.e.\
a surjective morphism of BCAs
$F((\ul{s}))[\sqbr{\ul{t}}] \ar A$.
Let $\ul{u} = (u_{1}, \ldots, u_{l})$ be a separating transcendency basis
for $F$ over $k$. By \cite{Ye1} Cor.\ 1.5.19,
$k \sqbr{\ul{u}, \ul{s}, \ul{t}} \ar F((\ul{s}))[\sqbr{\ul{t}}]$
is topologically \'{e}tale (rel.\ to $k$). Therefore, using \cite{Ye1}
formula (1.4.2) and Theorem 1.5.11, it follows that
\[ \cal{P}_{F((\ul{s}))[\sqbr{\ul{t}}]/k}^{n,\opn{sep}} \cong
F((\ul{s}))[\sqbr{\ul{t}}] \otimes_{k\sqbr{\ul{u}, \ul{s}, \ul{t}}}
\cal{P}_{k \sqbr{\ul{u}, \ul{s}, \ul{t}} / k}^{n,\opn{sep}} \]
is a free ST left $F((\ul{s}))[\sqbr{\ul{t}}]$-module of finite rank.
Now in general, if $\phi : M \surj N$ is a strict surjection of ST
$k$-modules, then so is
$\phi \otimes \phi : M \otimes_{k} M \ar N \otimes_{k} N$; and if
$M' \subset M$ and $N' \subset N$ and submodules such that
$\phi(M') \subset N'$, then
$\bar{\phi}: M / M' \surj N / N'$ is also strict.
This implies that
$\cal{P}_{F((\ul{s}))[\sqbr{\ul{t}}]/k}^{n,\opn{sep}} \ar
\cal{P}_{A/k}^{n,\opn{sep}}$ is a strict surjection. Hence
$\cal{P}_{A/k}^{n,\opn{sep}}$ is a ST module of finite type over
$A$ (as a left module, via $a \mapsto a \otimes 1$).

Given a finite type ST $A$-module $M$ as above, use Lemma \ref{lem4.3}.
\end{pf}


\section{$\protect\cal{D}$-Modules over TLFs}

Henceforth $k$ is a fixed perfect field.
Let $K$ be a topological local field(TLF) over $k$. We need to understand
the structure of the ring $\cal{D}(K)$ of continuous differential operators.
First assume $k$ has
characteristic $p$. Let $M$ be a free ST $K$-module of finite rank. We know
from \cite{Ye1} Theorems 2.1.14 and 1.4.9
that $\cal{D}(K; M)$ admits the $p$-filtration
\[ \cal{D}(K ;M) = \opn{Diff}_{K/k}(M,M) =
\bigcup_{n=0}^{\infty} \opn{End}_{K^{(p^{n}/k)}}(M)\ . \]
Here $K^{(p^{n}/k)} = k \otimes_{k} K$, with
$1 \otimes \lambda = \lambda^{p^{n}} \otimes 1$ for $\lambda \in k$.
This filtration is cofinal with the order filtration - see \cite{Ye1} Lemma
1.4.8. According to ibid.\ Prop.\ 2.1.13, the relative Frobenius map
$K^{(p^{n}/k)} \ar K$, $\lambda \otimes a \mapsto \lambda a^{p^{n}}$,
is a finite morphism in $\mathsf{TLF}(k)$.

In characteristic $0$, $\cal{D}(K)$ is a ``topologically \'{e}tale
localization'' of a Weyl algebra. Choose a parametrization
$K \cong F((\ul{s}))$ and a separating transcendency basis
$\ul{u}$ for $F$ over $k$. Let
$\ul{t} = (t_{1}, \ldots, t_{m}) := (\ul{u}, \ul{s}) =
(u_{1}, \ldots, s_{1} \ldots)$ be the concatenated sequence. Then
$k\sqbr{ \ul{t} } \ar K$ is a flat, topologically \'{e}tale homomorphism
in $\mathsf{STComAlg}(k)$. The ring
$\cal{D}(k \sqbr{\ul{t}})$ is a Weyl algebra over $k$:
$\cal{D}(k \sqbr{\ul{t}}) \cong k \sqbr{\ul{t}} \otimes_{k}
k \sqbr{ \partial_{1},\ldots,\partial_{m} }$, where
$\partial_{i}:= \frac{\partial}{\partial t_{i}}$, and the multiplication
is determined by
$(1 \otimes \partial_{i})(t_{j} \otimes 1)=t_{j}
\otimes \partial_{i} + (\partial_{i} * t_{j}) \otimes 1$. By Prop.\
\ref{prop4.1}, we have
$\cal{D}(K) \cong K \otimes_{k \sqbr{\ul{t}}}
\cal{D}(k \sqbr{\ul{t}})$. Considering the
faithful action of $\cal{D}(K)$ on $K$, we get a presentation
\begin{equation}  \label{eqn5.8}
\begin{array}{rcl}
\cal{D}(K) & \cong & K \otimes_{k} k \sqbr{\partial_{1},\ldots,\partial_{m}} \\
(1 \otimes \partial_{i})(a \otimes 1) & = &
a \otimes \partial_{i} + (\partial_{i} * a) \otimes 1
\end{array}
\end{equation}
for $i = 1, \ldots, m$ and $a \in K$
(i.e. $\cal{D}(K)$ is a smash product of $K$ and the universal enveloping
algebra of
the abelian $k$-Lie algebra spanned by the derivations $\partial_{i}$).

\begin{dfn}
Let $K$ be a TLF over $k$. Define
$\omega(K)$ to be the top degree component of
$\Omega^{\bdot, \mathrm{sep}}_{K/k}$.
It is a free ST $K$-module of rank $1$.
\end{dfn}

At this point we can exhibit the canonical right $\cal{D}(K)$-module
structure on $\omega(K)$ (cf.\ \cite{Bo} Ch.\ VI \S 3.2).

\begin{prop} \label{prop5.1}
For any  $K \in \mathsf{TLF}(k)$ there is a unique right
$\cal{D}(K)$-module structure on $\omega(K)$, written
$\alpha * D$, for $\alpha \in \omega(K)$ and $D \in \cal{D}(K)$,
such that:
\begin{enumerate}
\rmitem{i} If $D=a \in K$ then $\alpha * a = a \alpha$.
\rmitem{ii} If $D = \partial \in \cal{T}(K)$ then
$\alpha * \partial  = - \opn{L}_{\partial}(\alpha)$, where
$\opn{L}_{\partial}$ is the Lie derivative \rom{(}see Prop.\
\rom{\ref{prop4.2})}.
\rmitem{iii} If $\opn{char} k=p$ and
$D \in \cal{D}^{p^{n}-1}(K)$ for some $n \geq 0$, then for
every $a \in K$,
\[ \langle D * a, \alpha \rangle_{K/K^{(p^{n}/k)}} =
\langle a, \alpha * D \rangle_{K/K^{(p^{n}/k)}}\ , \]
where $\langle - , - \rangle_{K/K^{(p^{n}/k)}}$ is the trace pairing of
\cite{Ye1} formula \rom{(2.3.8)}.
\end{enumerate}
\end{prop}

\begin{pf}
First assume $\opn{char} k=0$. Since
$[\opn{L}_{\partial_{i}},\opn{L}_{\partial_{j}}]=0$,
$\alpha * \partial_{i} := - \opn{L}_{\partial_{i}}(\alpha)$ is an action
of
$k \sqbr{\partial_{1}, \ldots, \partial_{m}}$ on $\omega(K)$. According to
the presentation (\ref{eqn5.8}), in order to extend this to a
right action of $\cal{D}(K)$ it suffices to show that
\[ -a \opn{L}_{\partial_{i}}(\alpha) = - \opn{L}_{\partial_{i}}(a
\alpha) + \partial_{i}(a) \alpha \]
which is true since $\opn{L}_{\partial_{i}}$ is an even derivation of
$\Omega^{\bdot,\opn{sep}}_{K/k}$ and
$\opn{L}_{\partial_{i}}(a)= \partial_{i}(a)$. By Lemma \ref{lem4.1},
condition (ii) holds for an arbitrary
derivation $\partial = \sum a_{i} \partial_{i}$, $a_{i} \in K$.

Next consider the case $\opn{char} k=p$. Let
$D \in \cal{D}^{n}(K)$. By \cite{Ye1} Lemma 1.4.8, $D$ is
$K^{(p^{n}/k)}$-linear. The trace
pairing $\langle -,- \rangle_{K/K^{(p^{n}/k)}}$ is perfect (\cite{Ye1}
Prop.\ 2.3.9), so by
adjunction $D$ acts on $\omega(K)$. The functoriality of the trace
guarantees that this action is independent of $n$. We thus get a right
action satisfying conditions (i) and (iii). In order to check (ii) it
suffices to look at $\partial=\partial_{i}$.
Let $\beta:= \mathrm{d} t_{1} \wedge \cdots \wedge \mathrm{d} t_{i-1}
\wedge \mathrm{d} t_{i+1} \wedge \cdots \wedge \mathrm{d} t_{m}$.
We can compute the difference:
\begin{eqnarray*}
\lefteqn{\langle \partial_{i}(a),b \mathrm{d} t_{i} \wedge \beta
\rangle_{K/K^{(p^{n}/k)}}  -
\langle a,-\opn{L}_{\partial_{i}}(b \mathrm{d}
t_{i} \wedge \beta) \rangle_{K/K^{(p^{n}/k)}} } \blnk{30mm}\\
& & = \opn{Tr}_{K/K^{(p/k)}}(\partial_{i}(ab) \mathrm{d}
t_{i} \wedge \beta) \\
& & = \opn{Tr}_{K/K^{(p/k)}}(\mathrm{d}(ab \beta)) = 0
\end{eqnarray*}
since $\opn{Tr}_{K/K^{(p/k)}}$ commutes with $\mathrm{d}$ and vanishes
on $\Omega^{\leq m-1,\opn{sep}}_{K/k}$.
\end{pf}

Let $\cal{D}(K)^{\circ}$ denote the opposite ring of $\cal{D}(K)$.

\begin{prop} \label{prop5.11}
The right $\cal{D}(K)$ action on $\omega(K)$ of the previous proposition
induces a canonical isomorphism of filtered $k$-algebras
\[ T_{K} : \cal{D}(K)^{\circ} \iso
\cal{D}(K; \omega(K)). \]
\end{prop}

\begin{pf}
If $\opn{char} k = p$ we have, for every $n \geq 0$, an isomorphism (of
$K$-$K$-bimodules)
$T^{n}: \opn{End}_{K^{(p^{n}/k)}}(K)^{\circ}
\iso \opn{End}_{K^{(p^{n}/k)}}(\omega(K))$
induced by adjunction. In the limit we get $T_{K}$.

If $\opn{char} k = 0$, choose a topologically \'{e}tale homomorphism
$k \sqbr{\ul{t}} \ar K$. Let $T_{\ul{t}} : \cal{D}(K) \ar \cal{D}(K)$ be
the involution such that $T_{\ul{t}} |_{K}$ is the identity and
$T_{\ul{t}}(\partial_{i})=-\partial_{i}$ (cf.\ formula (\ref{eqn5.8})).
Let
$\phi_{\ul{t}}:K \iso \omega(K)$
be the $K$-linear isomorphism defined by
$\phi_{\ul{t}}(1)= \mathrm{d} t_{1} \wedge \cdots \wedge \mathrm{d} t_{m}$.
Then for any $D \in \cal{D}(K)$,
\[ T_{K}(D)= \phi_{\ul{t}} \circ T_{\ul{t}}(D) \circ \phi_{\ul{t}}^{-1}
\in \cal{D}(K; \omega(K)). \]
\end{pf}

\begin{cor} \label{cor5.10}
Let $K,\hat{K} \in \mathsf{BCA}(k)$ be fields and let $K \ar \hat{K}$ be a
topologically \'{e}tale  homomorphism in $\mathsf{STComAlg}(k)$.
Then $\omega(K) \ar \omega(\hat{K})$ is a homomorphism of right
$\cal{D}(K)$-modules.
\end{cor}

\begin{pf}
In characteristic $0$ this follows from the covariance of the Lie derivative.
In positive characteristics it follows from the fact that the trace map
commutes with base change, cf.\ \cite{Ye1} Prop.\ 2.3.11.
\end{pf}

On the category $\mathsf{TLF}(k)$ there is a functorial residue map. To
each morphism $f : K \ar L$ it assigns a homomorphism of differential
graded ST left $\Omega^{\bdot,\opn{sep}}_{K/k}$-modules,
$\opn{Res}_{L/K} = \opn{Res}_{f} : \Omega^{\bdot,\opn{sep}}_{L/k} \ar
\Omega^{\bdot,\opn{sep}}_{K/k}$ (cf.\ \cite{Ye1} Thm.\ 2.4.3).
The residue pairing
\[ \begin{array}{c}
\langle - , - \rangle_{L/K} : L \times \omega(L) \ar \omega(K) \\
\langle a , \alpha \rangle_{L/K} = \opn{Res}_{L/K}(a \alpha)
\end{array} \]
is a perfect pairing of ST $K$-modules, in the sense that the induced map
$\omega(L) \ar \opn{Hom}_{K}^{\opn{cont}}(L, \omega(K))$ is bijective
(cf.\ \cite{Ye1} Thm.\ 2.4.22 - Topological Duality).

\begin{thm} \label{thm5.1}
Let $K \in \mathsf{TLF}(k)$ and assume that $k \ar K$ is a morphism in
$\mathsf{TLF}(k)$. Given a DO  $D \in \cal{D}(K)$, let
$D^{\vee} \in \opn{End}_{k}(\omega(K))$
be its adjoint relative to the residue pairing
$\langle -,- \rangle_{K/k}$.
Then for every $\alpha \in \omega(K)$,
\[ D^{\vee}(\alpha) = \alpha * D . \]
In other words, the adjoint action of $\cal{D}(K)$ on $\omega(K)$
coincides with the canonical right action.
\end{thm}

\begin{pf}
We must show that for all $a \in K$, $\alpha \in \omega(K)$ and $D
\in \cal{D}(K)$, $\langle D * a, \alpha \rangle_{K/k} = \langle
a, \alpha * D \rangle_{K/k}$. In characteristic $p$ this follows
immediately from condition (iii) of Prop.\ \ref{prop5.1} and the
functoriality of the residue maps (\cite{Ye1} Thm.\ 2.4.2).

in characteristic $0$ first choose a parametrization
$K \cong F((\ul{t}))= F((t_{1},\ldots,t_{n}))$. Then $k \ar F$ is finite
separable
and any $k$-linear DO is also $F$-linear.
Given $\lambda \in F$,
$\ul{i} \in {\Bbb Z}^{n}$ and $1 \leq j \leq n$, write
$\opn{dlog}(\ul{t}) := \opn{dlog}(t_{1}) \wedge \cdots \wedge
\opn{dlog}(t_{n}) = t_{1}^{-1} \mathrm{d} t_{1} \wedge \cdots \wedge
t_{n}^{-1} \mathrm{d} t_{n}$
and
$\partial_{j} := \frac{\partial}{\partial t_{j}}$.
Then
\begin{eqnarray*}
\lefteqn{\opn{L}_{\partial_{j}}( \lambda \ul{t}^{\ul{i}}
\opn{dlog}(\ul{t}) ) =} \\
& & (-1)^{j-1} \mathrm{d} ( \lambda t_{j}^{-1} \ul{t}^{\ul{i}}
\opn{dlog}(t_{1}) \wedge \cdots \wedge \opn{dlog}(t_{j-1}) \wedge
\opn{dlog}(t_{j+1}) \wedge \cdots \wedge \opn{dlog}(t_{n}))
\end{eqnarray*}
so
$\opn{Res}_{K/k} ( \opn{L}_{\partial_{j}}( \lambda \ul{t}^{\ul{i}} )) =0$.
By continuity we conclude that
\begin{equation} \label{eqn5.10}
\opn{Res}_{K/k} ( \opn{L}_{\partial_{j}} (\alpha) ) =0
\end{equation}
for all $\alpha$.

To prove the theorem it suffices to consider either $D=b \in K$ or
$D=\partial_{j}$.
For $D=b$ we get
\[ \langle a, \alpha * b \rangle_{K/k} =
\langle a, b \alpha  \rangle_{K/k} = \opn{Res}_{K/k}(ab \alpha) =
\langle ab , \alpha \rangle_{K/k} =
\langle b * a, \alpha  \rangle_{K/k} . \]
For $D=\partial_{j}$ we use ``integration by parts'':
\begin{eqnarray*}
\lefteqn{ \langle \partial_{j} * a , \alpha \rangle_{K/k} -
 \langle  a , \alpha * \partial_{j} \rangle_{K/k} =} \\
& & \opn{Res}_{K/k}( \opn{L}_{\partial_{j}} (a) \alpha +
a \opn{L}_{\partial_{j}} (\alpha) ) =
\opn{Res}_{K/k}( \opn{L}_{\partial_{j}} (a \alpha)) = 0
\end{eqnarray*}
by (\ref{eqn5.10}).
\end{pf}


\section{Duals of Finite Type Modules}

The purpose of this subsection is to establish the existence of a canonical
dual module $\opn{Dual}_{A} M$ to every finite type ST $A$-module $M$.
If $k \ar A$ is a morphism in $\mathsf{BCA}(k)$, then we set
$\opn{Dual}_{A} M := \opn{Hom}_{k}^{\mathrm{cont}}(M, k)$, endowed with the
fine $A$-module topology. Otherwise we define $\opn{Dual}_{A} M$ using
differential operators, and show this definition is independent of choices
made by a base change argument, which reduces things to the case when
$k \ar A$ is a morphism.
Recall that $k$ is a fixed perfect field.
For a TLF $K$, $\omega(K)$ is the top degree component of
$\Omega^{\bdot, \mathrm{sep}}_{K/k}$, a rank $1$ free ST $K$-module.

\begin{dfn} \label{def6.1}
Let $A,K \in \mathsf{BCA}(k)$ be a local ring and a field,
respectively, and let $\sigma: K \ar A$ be a morphism in $\mathsf{BCA}(k)$.
For any finite type ST $A$-module $M$ define
\[ \opn{Dual}_{\sigma} M := \opn{Hom}_{K; \sigma}^{\opn{cont}}
(M, \omega(K)) , \]
the set of continuous $K$-linear homomorphisms, where $M$ is a $K$-module
via $\sigma$. Put on $\opn{Dual}_{\sigma} M$ the {\em fine $A$-module
topology}.
\end{dfn}

\begin{rem} \label{rem6.1}
The module $\opn{Hom}_{K; \sigma}^{\opn{cont}}(M, \omega(K))$, with the
(weak) $\opn{Hom}$ topology, is a ST $A$-module. Therefore the identity map
$\opn{Dual}_{\sigma} M \ar \opn{Hom}_{K; \sigma}^{\opn{cont}}(M, \omega(K))$
is continuous. However, this will not be a homeomorphism unless
$\sigma :K \ar A$ is a pseudo coefficient field and $M$ is a finite length
module.
\end{rem}

Let $A$ be a  commutative noetherian local ring, with maximal ideal $\mx$,
and let $I$ be an injective hull of $A / \mx$. Then
$M \mapsto \opn{Hom}_{A}(M,I)$ is a duality between finite type (i.e.\
finitely generated) $A$-modules and cofinite type (i.e.\ artinian)
$A$-modules. The module $\opn{Hom}_{A}(M,I)$ is called a {\em Matlis dual}
of $M$ (cf.\ \cite{LC} \S 4).

\begin{lem} \label{lem6.1}
Let $\sigma : K \ar A$ and $M$ be as in Def.\ \ref{def6.1}.
\begin{enumerate}
\rmitem{a} Suppose $\tau : L \ar A$ and $f : K \ar L$ are morphisms in
$\mathsf{BCA}(k)$, with $L$ a field, and $\sigma = \tau \circ f$. Then
the map
\[ \begin{array}{rcl}
\opn{Dual}_{\tau} M & \ar & \opn{Dual}_{\sigma} M \\
\phi & \mapsto & \opn{Res}_{L/K} \circ \phi
\end{array} \]
is an isomorphism of ST $A$-modules.

\rmitem{b} The (untopologized) $A$-module $\opn{Dual}_{\sigma} M$ is a Matlis
dual of $M$. In particular, Taking $M = A$, it follows that
$\opn{Dual}_{\sigma} A$ is an injective hull of $A / \mx$.
As a ST $A$-module, $\opn{Dual}_{\sigma} M$ is of cofinite type.
\end{enumerate}
\end{lem}

\begin{pf}
(a)\ First consider the case when $\tau: L \ar A$ is finite; so $A$ has the
fine $L$-module topology (Prop.\ \ref{prop2.1} (a)).
Then $M$ is a free ST $L$-module of
finite rank. By Topological Duality
(\cite{Ye1} Thm.\ 2.4.22),
$\opn{Dual}_{\tau} M \ar \opn{Dual}_{\sigma} M$
is bijective, and it is an isomorphism of ST $A$-modules since both modules
have the fine $A$-module topologies.

Next assume $\tau: L \ar A$ is a pseudo coefficient field. Because
$\omega(K)$ (resp.\ $\omega(L)$) is a simple, separated ST $K$-module
(resp.\ $L$-module), and
$M \cong \lim_{\leftarrow n} M / \mx^{n} M$,
we can use \cite{Ye1} Prop.\ 1.2.22 to conclude that
\begin{equation} \label{eqn6.5}
\opn{Dual}_{\sigma} M = \bigcup_{n = 1}^{\infty}
\opn{Hom}_{K; \sigma}^{\opn{cont}}(M / \mx^{n} M, \omega(K))
\end{equation}
and similarly for $L$.
For any $n \geq 1$, $M / \mx^{n} M$ is a finite type ST
$A / \mx^{n}$-module, so we are back to the first step.

For the general situation, we may factor $\tau$ through some pseudo
coefficient field $\tau' : L' \ar A$ (cf.\ Lemma \ref{lem2.2} (d)),
and use the functoriality of the residue maps.

\medskip \noindent (b)\
By part (a) we can assume that $\sigma : K \ar A$ is a pseudo coefficient
field. Then in (\ref{eqn6.5}) we can drop the superscript ``$\opn{cont}$'',
in which case the statement is well known (cf.\ \cite{LC} p.\ 63 Example 1).
\end{pf}

Let $A, \hat{A}$ be local BCAs, with maximal ideals
$\mx, \hat{\mx}$ respectively, and let $v : A \ar \hat{A}$ be an
intensification homomorphism. Note that $v$, being a local
homomorphism, is faithfully flat. Let $\sigma : K \ar A$ be a morphism
in $\mathsf{BCA}(k)$, with $K$ a field. Assume that there is an
intensification homomorphism $u : K \ar \hat{K}$
and a morphism $\hat{\sigma} : \hat{K} \ar \hat{A}$ s.t.\
$\hat{A} \cong A \otimes_{K}^{(\wedge)} \hat{K}$.

\begin{prop} \label{prop6.1}
Let $M$ be a finite type ST $A$-module, and set
$\hat{M} := \hat{A} \otimes_{A} M$. Then any
$\phi \in \opn{Dual}_{\sigma} M$ has a unique extension
$\hat{\phi} \in \opn{Dual}_{\hat{\sigma}} \hat{M}$. The resulting continuous
homomorphism
\[ q_{v; \sigma}^{M} : \opn{Dual}_{\sigma} M \ar
\opn{Dual}_{\hat{\sigma}} \hat{M} \]
is injective, and induces an isomorphism of ST $\hat{A}$-modules
\[ 1 \otimes q_{v; \sigma}^{M}:
\hat{A} \otimes_{A} \opn{Dual}_{\sigma} M \iso
\opn{Dual}_{\hat{\sigma}} \hat{M} . \]
\end{prop}

\begin{pf}
Let $n := \opn{res.dim} \sigma$. Then we can extend $\sigma$ to a
a pseudo coefficient field
$K((\ul{s})) = K((s_{1}, \ldots, s_{n})) \ar A$, and extend $\hat{\sigma}$
to $\hat{K}((\ul{s})) \ar \hat{A}$. By replacing $K, \hat{K}$ with
$K((\ul{s})), \hat{K}((\ul{s}))$ we can then assume that
$\sigma, \hat{\sigma}$  are pseudo coefficient fields. This puts us in the
setup of Prop.\ \ref{prop3.2}. For $i \geq 0$ define
\[ H^{i} := \opn{Hom}_{K; \sigma}^{\opn{cont}}(M / \mx^{i+1} M, \omega(K))
\subset \opn{Dual}_{\sigma} M\, \]
and similarly define $\hat{H}^{i}$.
Since $\opn{Dual}_{\sigma} M$ and $\opn{Dual}_{\hat{\sigma}} \hat{M}$ both have
the fine topologies, it suffices to exhibit an isomorphism
$\hat{A} \otimes_{A} H^{i} \iso \hat{H}^{i}$, with
$\hat{\phi} := 1 \otimes \phi$ extending $\phi$. By Prop.\ \ref{prop3.2} (a),
\[ \hat{K} \otimes_{K} (M / \mx^{i+1} M) \cong
\hat{A} \otimes_{A} (M / \mx^{i+1} M) \cong \hat{M} / \mx^{i+1} \hat{M} . \]
Since $K \ar \hat{K}$ is topologically \'{e}tale,
$\hat{K} \otimes_{K} \omega(K) \iso \omega(\hat{K})$. Therefore
$\hat{K} \otimes_{K} H^{i} \iso \hat{H}^{i}$; and again by Prop.\
\ref{prop3.2} (a), $\hat{A} \otimes_{A} H^{i} \iso \hat{H}^{i}$.
\end{pf}

Let $A$ be a local BCA with maximal ideal $\mx$. Suppose
$\sigma, \sigma' : K \ar A$ are pseudo coefficient fields, such that
$\sigma \equiv \sigma'\ (\opn{mod} \mx)$.
Let $M$ be a finite type ST $A$-module. Given a nonzero element
$x \in M$, its order with respect to $\mx$ is
\[ \opn{ord}_{\mx}(x) := \opn{max} \{ n\ |\ x \in \mx^{n} M \} . \]
If $\opn{ord}_{\mx}(x) =n$, then the symbol of $x$ is its image in
$\mx^{n} / \mx^{n+1} \subset \opn{gr}_{\mx} M$.

\begin{dfn} \label{def6.3}
An $\mx$-filtered $K$-basis of $M$ is
a sequence $\ul{x} = (x_{0}, x_{1}, \ldots)$ of elements of $M$,
such that the symbols of $x_{0}, x_{1}, \ldots$ form a $K$-basis of
$\opn{gr}_{\mx} M$, and such that
$\opn{ord}_{\frak{m}}(x_{i}) \leq \opn{ord}_{\frak{m}}(x_{i+1})$.
\end{dfn}

Choose such a basis $\ul{x}$.
Then any $x \in M$ is expressed uniquely as a convergent sum
\[ x = \sum_{i} \sigma'(\lambda_{i}) x_{i} =
\sum_{i} \sigma(\mu_{i}) x_{i} \]
with $\lambda_{i}, \mu_{i} \in K$. Define functions $D_{ij} : K \ar K$
by the equation
\[ \sigma'(\lambda) x_{i} = \sum_{j} \sigma(D_{ij}(\lambda)) x_{j} . \]

\begin{lem} \label{lem6.2}
$D_{ij} \in  \cal{D}(K)$, i.e.\ it is a continuous differential operator
over $K$ relative to $k$.
\end{lem}

\begin{pf}
Pick two indices $i_{0},i_{1}$, and let
$n:= \opn{max} \{ \opn{ord}_{\mx}(x_{i_{0}}), \opn{ord}_{\mx}(x_{i_{1}}) \}$.
We can  compute the function $D_{i_{0} i_{1}}$ for the module
$M / \mx^{n+1} M$ instead of $M$. Define
\[ A^{-} := \sigma(K) \oplus \mx = \sigma'(K) \oplus \mx \subset A . \]
This is a local BCA, with $A^{-} / \mx \cong K$,
and $A^{-} \ar A$ is a finite morphism. Let $l$ be the length of
$M / \mx^{n+1} M$ over $A^{-}$, and let
$E,E' : K^{l} \iso M / \mx^{n+1} M$ be the $K$-linear homeomorphisms
\begin{eqnarray*}
E(\lambda_{0}, \ldots, \lambda_{l-1}) & := &
\sum_{i=0}^{l-1} \sigma(\lambda_{i}) x_{i} \\
E'(\lambda_{0}, \ldots, \lambda_{l-1}) & := &
\sum_{i=0}^{l-1} \sigma'(\lambda_{i}) x_{i} \\
\end{eqnarray*}
($M / \mx^{n+1} M$ is a free ST $K$-module via $\sigma$ and via $\sigma'$).
According to \cite{Ye1} Prop.\ 1.4.4, $E,E^{-1},E'$ and $(E')^{-1}$ are
DOs over $A^{-}$, relative to $k$. Set
\[ D:= E^{-1} \circ E' : K^{l} \iso K^{l} \ , \]
which is a DO over $A^{-}$, and hence over $K$. Expanding $D$ as an
$l \times l$ matrix with entries in $\cal{D}(K)$, one gets
$D = [D_{ij}]$.
\end{pf}

One can easily show that
\[ \opn{ord}_{K}(D_{ij}) \leq 2 \opn{max} \{ -1,
\opn{ord}_{\frak{m}}(x_{j}) - \opn{ord}_{\frak{m}}(x_{i}) \} \]
and
$D_{ii} = 1$. Thus the matrix of DOs looks like this:
\[ [ D_{ij} ] =
\left[ \begin{array}{cccc}
1 & * & * \\
0 & 1 & * & \cdot \ \cdot \\
0 & 0 & 1 \\
  & : \\
\end{array} \right] \]

\begin{dfn} \label{dfn6.2}
In the situation described above, define a function
\[ \Psi_{\sigma, \sigma'}^{M} : \opn{Dual}_{\sigma} M \ar \opn{Dual}_{\sigma'}
M \]
by the equation
\begin{equation} \label{eqn6.6}
\Psi_{\sigma, \sigma'}^{M}(\phi) (\sum_{i} \sigma'(\lambda_{i}) x_{i}) =
\sum_{i,j} \lambda_{i} (\phi(x_{j}) * D_{ij})
\end{equation}
for $\phi \in \opn{Dual}_{\sigma} M$ and $\lambda_{i} \in K$.
\end{dfn}

The second sum in (\ref{eqn6.6}) makes sense, since there are only finitely
many nonzero terms in it. At first glance this somewhat strange definition
seems to depend on the basis $\ul{x}$. We shall soon see that there is no
dependence on the basis, and that in fact $\Psi_{\sigma, \sigma'}^{M}$ is an
isomorphism of ST $A$-modules.
Immediately from the definition we get:

\begin{lem} \label{lem6.0}
$\Psi^{M}_{\sigma, \sigma'}$ is a $k$-linear bijection, with inverse
$\Psi^{M}_{\sigma', \sigma}$. Given a third pseudo coefficient field
$\sigma'': K \ar A$ one has
\[ \Psi^{M}_{\sigma, \sigma''} =
\Psi^{M}_{\sigma'', \sigma'} \circ \Psi^{M}_{\sigma, \sigma'}. \]
\end{lem}

Further properties of  $\Psi^{M}_{\sigma, \sigma'}$ are less obvious.

\begin{lem} \label{lem6.3}
Under the combined assumptions of Prop.\ \ref{prop6.1} and Def.\
\ref{dfn6.2}, one has
\[ \Psi_{\hat{\sigma}, \hat{\sigma}'}^{\hat{M}} \circ
q_{v; \sigma}^{M} = q_{v; \sigma'}^{M} \circ
\Psi_{\sigma, \sigma'}^{M}. \]
Here we are using the $\hat{\mx}$-filtered basis
$(1 \otimes x_{0}, 1 \otimes x_{1}, \ldots)$ on $\hat{M}$ to define
$\Psi_{\hat{\sigma}, \hat{\sigma}'}^{\hat{M}}$.
\end{lem}

\begin{pf}
The DOs $\hat{D}_{ij} \in \cal{D}(\hat{K})$ which appear in the definition of
$\Psi_{\hat{\sigma}, \hat{\sigma}'}^{\hat{M}}$ are precisely the images
of the DOs
$D_{ij} \in \cal{D}(K)$ under the natural ring homomorphism
$\cal{D}(K) \ar \cal{D}(\hat{K})$. By Cor.\ \ref{cor5.10},
$\omega(K) \ar \omega(\hat{K})$ is a homomorphism of right
$\cal{D}(K)$-modules.
\end{pf}

\begin{lem} \label{lem6.4}
In the situation of Def.\ \ref{dfn6.2},
suppose in addition that $k \ar K$ is a morphism in $\mathsf{BCA}(k)$. Then
for any $\phi \in \opn{Dual}_{\sigma} M$, one has
\[ \opn{Res}_{K/k} \circ \phi =
\opn{Res}_{K/k} \circ \Psi_{\sigma, \sigma'}^{M}(\phi). \]
\end{lem}

\begin{pf}
Say $\phi(x_{i}) = \alpha_{i} \in \omega(K)$. Given
$x = \sum_{i} \sigma'(\lambda_{i}) x_{i} \in M$, with
$\lambda_{i} \in K$, we have by definition
$x = \sum_{i,j} \sigma(D_{ij} * \lambda_{i}) x_{j}$. So
\[ \opn{Res}_{K/k} \circ \phi (x) =
\opn{Res}_{K/k} (\sum_{i,j} (D_{ij} * \lambda_{i}) \alpha_{j})  . \]
On the other hand, setting
$\phi' := \Psi_{\sigma, \sigma'}^{M}(\phi)$, one has
\[ \opn{Res}_{K/k} \circ \phi' (x) =
\opn{Res}_{K/k} (\sum_{i,j} \lambda_{i} (\alpha_{j} * D_{ij}))  . \]
By linearity and continuity, it suffices to prove that for all
$i,j \geq 0$:
\[ \opn{Res}_{K/k} \left( (D_{ij} * \lambda_{i}) \alpha_{j} \right)
= \opn{Res}_{K/k} \left( \lambda_{i} (\alpha_{j} * D_{ij}) \right)\ ; \]
but this is done in Thm.\ \ref{thm5.1}.
\end{pf}

\begin{lem} \label{lem6.5}
Let $K \in \mathsf{BCA}(k)$ be a field. There exists a field
$\hat{K} \in \mathsf{BCA}(k)$ and a homomorphism $u : K \ar \hat{K}$
in $\mathsf{STComAlg}(k)$ such that $k \ar \hat{K}$ is a morphism
in $\mathsf{BCA}(k)$ and $u$ is an intensification. Moreover we can choose
$u$ to be dense.
\end{lem}

\begin{pf}
Choose a parametrization $K \cong F((\ul{t}))$. $F$ is a finitely generated
field extension of $k$; let $\ul{s} = (s_{1}, \ldots, s_{m})$ be a
transcendency basis for $F/k$. Then $k(\ul{s}) \ar F$ is a finite morphism
in $\mathsf{BCA}(k)$. The map
$k(\ul{s}) \ar k((\ul{s}))$ is certainly an intensification homomorphism.
Applying finitely ramified base change (Thm.\ \ref{thm3.1})
we get a dense intensification homomorphism
$K \ar K \otimes_{k(\ul{s})}^{(\wedge)} k((\ul{s}))$. Thus the BCA
$K \otimes_{k(\ul{s})}^{(\wedge)} k((\ul{s}))$ is a reduced cluster of TLFs,
and we can take $\hat{K}$ to be any local factor of it.
\end{pf}

\begin{prop} \label{prop6.2}
Let $A$ be a local BCA with maximal ideal $\mx$, and let
$\sigma, \sigma' : K \ar A$ be two pseudo coefficient fields, such that
$\sigma \equiv \sigma'\ (\opn{mod} \mx)$. Let $M$ be a finite type ST
$A$-module. Then the map
$\Psi_{\sigma, \sigma'}^{M}$
is an isomorphism of ST $A$-modules, independent of the $\mx$-filtered
$K$-basis
$\ul{x} = (x_{0}, x_{1}, \ldots)$.
\end{prop}

\begin{pf}
First we reduce the problem to the case when $K \iso A / \mx$, i.e.\
when $\sigma, \sigma'$ are coefficient fields. Let $A^{-}$ be the algebra
$\sigma(K) \oplus \mx \subset A$, cf.\ proof of Lemma \ref{lem6.2}.
The map $\Psi_{\sigma, \sigma'}^{M}$ is the
same when restricting $M$ to an $A^{-}$-module, so we may replace $A$ with
$A^{-}$.

Choose an intensification homomorphism $u : K \ar \hat{K}$ as in Lemma
\ref{lem6.5}, and define $\hat{A} := A \otimes_{K}^{(\wedge)} \hat{K}$,
w.r.t.\ the morphism $\sigma : K \ar A$. So the homomorphism
$v : A \ar \hat{A}$ is also an intensification,
$\hat{A}$ is local with maximal ideal
$\hat{\mx} = \hat{A} \cdot v(\mx)$, and
$\hat{A} / \hat{\mx} \cong \hat{K}$.  Let
$R: \opn{Dual}_{\hat{\sigma}} \hat{M} \iso
\opn{Hom}_{k}^{\mathrm{cont}}(\hat{M}, k)$
be the $\hat{A}$-linear isomorphism
$\phi \mapsto \opn{Res}_{\hat{K} / k} \circ \phi$ of Lemma \ref{lem6.1},
and similarly define $R'$.
According to  Lemmas \ref{lem6.3} and \ref{lem6.4}, the diagram

\medskip \noindent
\begin{equation} \label{eqn6.10}
\setlength{\unitlength}{0.30mm}
\begin{array}{ccccc}
\opn{Dual}_{\sigma} M & \lrar{q_{v; \sigma}^{M}}
& \opn{Dual}_{\hat{\sigma}} \hat{M} \\
\ldar{\Psi_{\sigma, \sigma'}^{M}} &
& \ldar{\Psi_{\hat{\sigma}, \hat{\sigma}'}^{\hat{M}}}
& \ldrar{R} \\
\opn{Dual}_{\sigma'} M & \lrar{q_{v; \sigma'}^{M}}
& \opn{Dual}_{\hat{\sigma}'} \hat{M} & \lrar{R'}
& \opn{Hom}_{k}^{\mathrm{cont}}(M, k)
\end{array}
\end{equation}
is commutative. Since $q_{v; \sigma}^{M}$ and $q_{v; \sigma'}^{M}$
are injections, we deduce
the independence of $\Psi_{\sigma, \sigma'}^{M}$ of the basis
$\ul{x}$, and that $\Psi_{\sigma, \sigma'}^{M}$
is an $A$-linear bijection. Since both
$\opn{Dual}_{\sigma} M$ and $\opn{Dual}_{\sigma'} M$ have the fine topologies,
$\Psi_{\sigma, \sigma'}^{M}$ is in fact a homeomorphism.
\end{pf}

\begin{prop} \label{prop6.3}
Under the hypothesis of Prop.\ \ref{prop6.2},
suppose $\tau, \tau' : L \ar A$ are pseudo coefficient fields,
and $f : K \ar L$ is a (finite) morphism in $\mathsf{BCA}(k)$, such that
$\tau \equiv \tau'\ (\opn{mod} \mx)$, $\sigma = \tau \circ f$ and
$\sigma' = \tau' \circ f$. Then for any $\phi \in \opn{Dual}_{\tau} M$ one has
\begin{equation} \label{eqn6.7}
\Psi_{\sigma, \sigma'}^{M} (\opn{Tr}_{f} \circ \phi) =
\opn{Tr}_{f} \circ \Psi_{\tau, \tau'}^{M}(\phi) .
\end{equation}
\end{prop}

\begin{pf}
After making a reduction as in prop.\ \ref{prop6.2}, we can assume that
$L \iso A / \mx$. Now set
$A^{-} := \sigma(K) \oplus \mx \subset A$. Choose a homomorphism
$u : K \ar \hat{K}$ as in Lemma \ref{lem6.5}, and define BCAs
$\hat{A}^{-} := A^{-} \otimes_{K}^{(\wedge)} \hat{K}$ and
$\hat{A} := A \otimes_{K}^{(\wedge)} \hat{K}$, w.r.t.\ the morphisms
$\sigma : K \ar A^{-} \ar A$. Let $v : A \ar \hat{A}$ be the resulting
intensification homomorphism. The algebra $\hat{A}^{-}$ is local, with
maximal ideal $\hat{A}^{-} \cdot v(\mx)$. Denote by $\hat{\frak{r}}$ the
Jacobson radical of $\hat{A}$; so
$\hat{\frak{r}} = \hat{A} \cdot v(\mx)$. Set
$\hat{L} := \hat{A} / \hat{\frak{r}} \cong L \otimes_{K} \hat{K}$.
For each $\hat{\mx} \in \opn{Max} \hat{A}$ denote by
$f_{\hat{\mx}} : K \ar L_{\hat{\mx}}$,
$v_{\hat{\mx}} : A \ar \hat{A}_{\hat{\mx}}$ and
$u_{\hat{\mx}} : L \ar \hat{L}_{\hat{\mx}}$ the localized homomorphisms.
We have
$\hat{A}_{\hat{\mx}} \cong A \otimes_{L}^{(\wedge)} \hat{L}_{\hat{\mx}}$,
and there are coefficient fields
$\hat{\tau}_{\hat{\mx}}, \hat{\tau}'_{\hat{\mx}} :
\hat{L}_{\hat{\mx}} \ar \hat{A}_{\hat{\mx}}$
extending $\tau, \tau'$. All the claims above follow from Prop.\
\ref{prop3.2}.

Let $\hat{M} := \hat{A} \otimes_{A} M$. For every
$\hat{\mx} \in \opn{Max} \hat{A}$ there is a homomorphism
\[ q_{v_{\hat{\mx}}; \tau}^{M} : \opn{Dual}_{\tau} M \ar
\opn{Dual}_{\hat{\tau}_{\hat{\mx}}} \hat{M}_{\hat{\mx}}\, \]
and a corresponding homomorphism $q_{v_{\hat{\mx}}; \tau'}^{M}$,
which, by Lemma \ref{lem6.3}, intertwine
$\Psi_{\tau, \tau'}^{M}$ with
$\Psi_{\hat{\tau}_{\hat{\mx}}, \hat{\tau}_{\hat{\mx}}'}^{
\hat{M}_{\hat{\mx}}}$.
There are also (injective) homomorphisms
$q_{v; \sigma}^{M}$ and $q_{v; \sigma'}^{M}$.
Since the trace maps satisfy
\[ u \circ \opn{Tr}_{f} = \sum_{\hat{\mx}}
\opn{Tr}_{f_{\hat{\mx}}} \circ u_{\hat{\mx}} \]
we get
\[ q_{v; \sigma}^{M} (\opn{Tr}_{f} \circ \phi) =
\sum_{\hat{\mx}} \opn{Tr}_{f_{\hat{\mx}}} \circ q_{v_{\hat{\mx}}; \tau}^{M}
(\phi) \]
and similarly with $\sigma', \tau'$,
so the problem is reduced to the case when $k \ar K$ is a morphism.

In this case, using \ref{lem6.4}  twice and the transitivity of residues,
we get
\[ \opn{Res}_{K/k} \circ \Psi_{\sigma, \sigma'}^{M} (\opn{Tr}_{f} \circ \phi)
= \opn{Res}_{K/k} \circ \opn{Tr}_{f} \circ \Psi_{\tau, \tau'}^{M}(\phi) \]
which, in virtue of Lemma \ref{lem6.1} (a), implies formula
(\ref{eqn6.7}).
\end{pf}

We are ready to prove the first main result of this article.

\begin{thm} \label{thm6.1} \rom{(Dual Modules)}\
Let $A$ be a local Beilinson completion algebra over $k$, and let $M$ be a
finite type semi topological $A$-module.
Then the following data exist:
\begin{enumerate}
\rmitem{a} A ST $A$-module
$\opn{Dual}_{A} M$, called the {\em dual module} of $M$.

\rmitem{b} For every morphism
$\sigma : K \ar A$ in $\mathsf{BCA}(k)$, with $K$ a field, an
isomorphism of ST $A$-modules
\[ \Psi_{\sigma}^{M} : \opn{Dual}_{A} M \iso \opn{Dual}_{\sigma} M =
\opn{Hom}^{\mathrm{cont}}_{K; \sigma}(M, \omega(K)). \]
\end{enumerate}

These data satisfy, and are completely determined by the following
conditions:

\begin{enumerate}
\rmitem{i} Let $f : K \ar L$ and $\tau : L \ar A$ be
morphisms in $\mathsf{BCA}(k)$, with $K,L$ fields, and let
$\sigma := \tau \circ f$. Then for any $\phi \in \opn{Dual}_{A} M$,
\[ \Psi_{\sigma}^{M}(\phi) = \opn{Res}_{f} \circ \Psi_{\tau}^{M}(\phi). \]
Here
$\opn{Res}_{f}: \omega(L) \ar \omega(K)$ is the residue map in
$\mathsf{TLF}(k)$, cf.\ \cite{Ye1} \S \rom{2.4}.
\rmitem{ii} Denote by $\mx$ the maximal ideal of $A$.
If $\sigma,\sigma' : K \ar A$ are pseudo coefficient fields such that
$\sigma \equiv \sigma'\ (\opn{mod} \mx)$, then
\[ \Psi_{\sigma'}^{M} = \Psi_{\sigma,\sigma'}^{M} \circ \Psi_{\sigma}^{M}, \]
where $\Psi_{\sigma,\sigma'}^{M}$ is the isomorphism defined in
Def.\ \ref{dfn6.2}.
\end{enumerate}
\end{thm}

Observe that if $A=K$ is a TLF, then there
is a canonical isomorphism
$\opn{Dual}_{K} M$ \linebreak
$\cong \opn{Hom}_{K}(M, \omega(K))$,
corresponding to the identity morphism $K \ar K$, thought of as a
coefficient field.

\begin{pf}
The proof is divided into four steps.

\medskip \noindent (1)\
Fix a coefficient field
$\tau_{0} : L_{0} = A / \mx \ar A$, and set
$\opn{Dual} M := \opn{Dual}_{\tau_{0}} M$. Given another coefficient field
$\tau : L_{0} \ar A$, we are forced by condition (ii) to define
$\Psi_{\tau}^{M} := \Psi_{\tau_{0}, \tau}^{M}$. For any other coefficient field
$\tau' : L_{0} \ar A$ this condition is satisfied, on account of Lemma
\ref{lem6.0}; condition (i) is irrelevant.

\medskip \noindent (2)\
Now let $\sigma : K \ar A$ be a pseudo coefficient field which factors
through some coefficient field
$\tau : L_{0} \ar A$ (if $\sigma$ is a quasi coefficient field then there is
precisely one such $\tau$). Define
$\Psi_{\sigma}^{M} : \opn{Dual} M \iso \opn{Dual}_{\sigma}$ to be
$\phi \mapsto \opn{Tr}_{L_{0} / K} \circ \Psi_{\tau_{0}, \tau}^{M}(\phi)$,
as is forced by condition (i). According to Prop.\ \ref{prop6.3}, this
definition is independent of the coefficient field $\tau$.

\medskip \noindent (3)\
Let $\sigma : K \ar A$ be any pseudo coefficient field. Choose some
pseudo coefficient field $\sigma' : K \ar A$ such that
$\sigma \equiv \sigma'\ (\opn{mod} \mx)$ and such that $\sigma'$ factors
through some coefficient field. For example, take
$\sigma' := \tau_{0} \circ \pi \circ \sigma$, where $\pi : A \surj L_{0}$
is the natural projection. Define
$\Psi_{\sigma}^{M} := \Psi_{\sigma', \sigma}^{M} \circ
\Psi_{\sigma'}^{M}$. Prop.\ \ref{prop6.3} shows that this definition is
independent of the choice of $\sigma'$, and furthermore it shows that
conditions (i) and (ii) hold for all pseudo coefficient fields.

\medskip \noindent (4)\
Finally let $\sigma : K \ar A$ be a morphism with
$\opn{res.dim} \sigma \geq 1$. Choose a
factorization $\sigma = \tau \circ f$, with
$\tau : L \ar A$ a pseudo coefficient field and $f : K \ar L$ a morphism.
Define
$\Psi_{\sigma}^{M}(\phi) :=  \opn{Res}_{f} \circ \Psi_{\tau}^{M}(\phi)$,
$\phi \in \opn{Dual} M$. Now condition (ii) is no longer relevant.
To verify condition (i) it suffices to prove the
independence of this definition on $\tau$. So suppose that $\sigma$ also
factors into $\sigma = \tau' \circ f'$.

First assume there exists some finite morphism $g : L \ar L'$ such that
$\tau = \tau' \circ g$ and $f' = g \circ f$. Then applying condition (i) to
$\tau = \tau' \circ g$, we get
\begin{equation} \label{eqn6.9}
\opn{Res}_{f} \circ \Psi_{\tau}^{M}(\phi) =
\opn{Res}_{f} \circ \opn{Tr}_{g} \circ \Psi_{\tau'}^{M}(\phi) =
\opn{Res}_{f'} \circ \Psi_{\tau'}^{M}(\phi)
\end{equation}
for $\phi \in \opn{Dual} M$. By taking $L'$ to be the separable closure of
$L$ in $L_{0}$, and then using formula (\ref{eqn6.9}), we can assume that
$L \ar L_{0}$ is purely inseparable.

It remains to consider the case when $L,L' \subset L_{0}$, and both
$L \ar L_{0}$ and $L' \ar L_{0}$ are purely inseparable. Choose $j >> 0$
such that
$L_{0}^{(p^{j}/k)} \subset L \cap L'$. Define
$L_{1} := K L_{0}^{(p^{j}/k)} \subset L_{0}$ and let
$\tau_{1}, \tau_{1}' : L_{1} \ar A$
be the restrictions of $\tau, \tau'$.  To finish the verification use formula
(\ref{eqn6.9}) twice more.
\end{pf}


\section{Traces on Dual Modules}

As before, $k$ is a fixed perfect field.
Suppose $A$ is a local BCA.
Then $\opn{Dual}_{A}: M \mapsto \opn{Dual}_{A} M$ is a functor on the category
of finite type ST $A$-modules.
Given a finite type ST $A$-module $M$ and an element $x \in M$, let
$\rho_{x} : A \ar M$ be the function $a \mapsto ax$. As in \cite{LC} Lemma
4.1, and by our Lemma \ref{lem6.1} (b),
sending $\phi \in \opn{Dual}_{A} M$ to the homomorphism
$x \mapsto \opn{Dual}_{A}(\rho_{x})(\phi)$ gives an isomorphism
$\opn{Dual}_{A} M \ar \opn{Hom}_{A}(M, \cal{K}(A))$.

Any BCA $A$ over $k$ decomposes into local factors:
$A = \prod_{\mx \in \opn{Max} A} A_{\mx}$, as ST $k$-algebras. Any
morphism in $\mathsf{BCA}(k)$ decomposes accordingly.

\begin{dfn} \label{dfn7.3}
Let $A$ be a BCA over $k$. Define
\[ \cal{K}(A) := \bigoplus_{\mx \in \opn{Max} A} \opn{Dual}_{A_{\mx}} A_{\mx}.
\]
Given any ST $A$-module $M$, define
\[ \opn{Dual}_{A} M := \opn{Hom}_{A}^{\opn{cont}}(M, \cal{K}(A)) \]
with the (weak) $\opn{Hom}$ topology.
\end{dfn}

With this definition $\opn{Dual}_{A}$ is an additive functor
$\mathsf{STMod}(A)^{\circ} \ar \mathsf{STMod}(A)$.
In view of the previous discussion and Prop.\ \ref{prop1.1} (2), there is
no conflict of definitions when
$A$ is local and $M$ is a ST $A$-module of finite type.

\begin{prop} \label{prop7.1} \rom{(Covariance of dual modules)}\
Let $v: A \ar \hat{A}$ be an intensification homomorphism between two BCAs.
Given a ST $A$-module $M$, set $\hat{M} := \hat{A} \otimes_{A} M$. Then
there is a unique homomorphism in $\mathsf{STMod}(A)$,
\[ q_{v}^{M} : \opn{Dual}_{A} M \ar  \opn{Dual}_{\hat{A}} \hat{M}, \]
with the following properties:
\begin{enumerate}
\rmitem{i} If $\phi: M \ar N$ is a homomorphism in $\mathsf{STMod}(A)$,
then
\[ q_{v}^{M} \circ \opn{Dual}_{A}(\phi) =
\opn{Dual}_{\hat{A}} (1 \otimes \phi) \circ q^{N}_{v}. \]
In other words,
$q_{v}: \opn{Dual}_{A} \ar \opn{Dual}_{\hat{A}} (\hat{A} \otimes_{A} -)$
is a natural transformation of functors.

\rmitem{ii} If $M$ is a ST $A$-module of finite type then the induced
homomorphism
\[ 1 \otimes q^{M}_{v} : \hat{A} \otimes_{A} \opn{Dual}_{A} M
\ar \opn{Dual}_{\hat{A}} \hat{M} \]
is an isomorphism.

\rmitem{iii} Let $\sigma : K \ar A$ be a morphism in $\mathsf{BCA}(k)$.
Suppose $K$ is a field, $A$ is local, and there is an intensification
homomorphism $u : K \ar \hat{K}$ s.t.\
$\hat{A} \cong A \otimes_{K}^{(\wedge)} \hat{K}$. Then for any ST $A$-module
of finite type $M$,
\[ q^{M}_{v} = (\Psi_{\hat{\sigma}}^{\hat{M}})^{-1} \circ
q^{M}_{v; \sigma} \circ \Psi_{\sigma}^{M} . \]

\rmitem{iv} If $w: \hat{A} \ar \Hat{\Hat{A}}$ is another
flat, finitely ramified, radically unramified and topologically
\'{e}tale homomorphism, then
$q_{w \circ v} = q_{w} \circ q_{v}$.
\end{enumerate}
These properties characterize $q^{M}_{v}$.
\end{prop}

\begin{pf}
We may assume $A$ is local. Let us first check uniqueness. If $M$ is a
finite type ST $A$-module, this follows from condition (iii). If $M$ has
the fine topology then $M \cong \lim_{\alpha \ar} M_{\alpha}$ with each
$M_{\alpha}$ a finite type module. By Lemma \ref{lem1.1} (4) we get
$\opn{Dual}_{A} M \cong \lim_{\leftarrow \alpha} \opn{Dual}_{A} M_{\alpha}$,
and we may use condition (i). Finally any ST $A$-module $M$ is a quotient of a
module $\tilde{M}$ which has the fine topology, and
$\opn{Dual}_{A} M \inj \opn{Dual}_{A} \tilde{M}$.

To define $q^{M}_{v}$ for $M$ of finite type amounts, essentially,
to repeating the steps of the proof of Thm.\ \ref{thm6.1}, using Lemma
\ref{lem6.3} at every step. For a general ST $A$-module $M$, let
$q^{M}_{v}$ be the canonical continuous homomorphism
\[ \opn{Hom}_{A}^{\opn{cont}}(M, \cal{K}(A)) \ar
\opn{Hom}_{\hat{A}}^{\opn{cont}}(\hat{A} \otimes_{A} M, \cal{K}(\hat{A})) \]
induced by
$q_{v} = q_{v}^{A} : \cal{K}(A) \ar \cal{K}(\hat{A})$.
\end{pf}

\begin{dfn} \label{dfn7.2}
Let $K,A \in \mathsf{BCA}(k)$, with $K$ a field, and let
$\sigma : K \ar A$ be a morphism. Define
\[ \opn{Res}_{A/K} = \opn{Res}_{\sigma} : \cal{K}(A) \ar
\omega(K) \]
to be the function sending
\[ \phi = \sum_{\frak{m}} \phi_{\frak{m}} \in \cal{K}(A) =
\bigoplus_{\frak{m}} \opn{Dual}_{A_{\frak{m}}} A_{\frak{m}} \]
to
$\sum_{\frak{m}} \Psi_{\sigma}^{A_{\frak{m}}}(\phi_{\frak{m}})(1) \in
\omega(K)$. Here $\frak{m}$ runs through the maximal ideals of $A$.
\end{dfn}

The residue map $\opn{Res}_{A/K}$ is $K$-linear. It is also continuous:
this follows from the adjunction formula, Lemma \ref{lem1.2}
(cf.\ Remark \ref{rem6.1}).
Because of the transitivity of residues, if there is a factorization
$\sigma: K \exar{ f } L \exar{ \tau } A$, then
$\opn{Res}_{A/K}= \opn{Res}_{L/K} \circ \opn{Res}_{A/L}$.

Here is the second main result of this article:

\begin{thm} \label{thm7.2}
\rom{(Traces)}\ Let $f : A \ar B$ be a morphism in $\mathsf{BCA}(k)$.
There is a unique continuous $A$-linear homomorphism
\[ \opn{Tr}_{B/A} = \opn{Tr}_{f} : \cal{K}(B) \ar \cal{K}(A) \]
having the following properties:

\begin{enumerate}
\rmitem{i} (Transitivity)\ Given another morphism $g : B \ar C$, one has
\[ \opn{Tr}_{C/A} = \opn{Tr}_{B/A} \circ \opn{Tr}_{C/B} . \]

\rmitem{ii} (Base Change )\ Suppose $u : A \ar \hat{A}$ is an intensification
homomorphism. Let
$\hat{B} := B \otimes_{A}^{(\wedge)} \hat{A}$, $v : B \ar \hat{B}$
and $\hat{f} : \hat{A} \ar \hat{B}$ be the algebras and homomorphisms gotten
by intensification base change (cf.\ Thm.\ \ref{thm3.1}). Then
\[ q_{u} \circ \opn{Tr}_{B/A} =
\opn{Tr}_{\hat{B}/\hat{A}} \circ q_{v}, \]
where $q_{u}, q_{v}$ are the homomorphisms of Prop.\ \ref{prop7.1}.

\rmitem{iii} If $A$ is a field, then
$\opn{Tr}_{B/A} = \opn{Res}_{B/A} : \cal{K}(B) \ar \cal{K}(A) = \omega(A)$.

\rmitem{iv} The map
\[ \cal{K}(B) \ar \opn{Hom}_{A}^{\mathrm{cont}}(B, \cal{K}(A)) \]
induced by $\opn{Tr}_{B/A}$ is bijective.
\end{enumerate}
\end{thm}

\begin{pf}
We may assume both $A,B$ are local, with maximal ideals $\mx, \frak{n}$.
Given any morphism $\sigma: K \ar A$ with $K$ a field, define
$\opn{Tr}_{B/A; \sigma} : \cal{K}(B) \ar \cal{K}(A)$
by
\begin{equation} \label{eqn7.2}
\opn{Tr}_{B/A; \sigma} := (\Psi_{\sigma}^{A})^{-1} \circ
\opn{Dual}_{\sigma}(f) \circ \Psi_{f \circ \sigma}^{B}
\end{equation}
where $\opn{Dual}_{\sigma}(f)(\phi) = \phi \circ f$ for
$\phi \in \opn{Dual}_{f \circ \sigma} B$.

The claim is that $\opn{Tr}_{B/A; \sigma}$ is independent of $\sigma$.
Let $\tau : L = A / \mx \ar A$ be any coefficient field. It suffices to
prove that
$\opn{Tr}_{B/A; \sigma} = \opn{Tr}_{B/A; \tau}$. To do so we choose an
intensification homomorphism $K \ar \hat{K}$ s.t.\
$k \ar \hat{K}$ is a morphism of BCAs, and set
$\hat{A} := A \otimes_{K}^{(\wedge)} \hat{K}$,
$\hat{B} := B \otimes_{K}^{(\wedge)} \hat{K}$ and
$\hat{L} := L \otimes_{K}^{(\wedge)} \hat{K}$.
Let $\hat{\tau} : \hat{L} \ar \hat{A}$ be the unique extension of $\tau$.
Note that by Prop.\ \ref{prop3.2},
$\hat{A} \cong A \otimes_{L}^{(\wedge)} \hat{L}$.
According to Prop.\ \ref{prop7.1} (iii),
\[ q_{\hat{A} / A} \circ \opn{Tr}_{B / A; \sigma} =
\opn{Tr}_{\hat{B} / \hat{A}; \hat{\sigma}} \circ q_{\hat{B} / B}  \]
and similarly for $\tau$. Since $\rho : k \ar \hat{A}$ is a morphism,
we get (using Thm.\ \ref{thm6.1} (i))
\[ \opn{Tr}_{\hat{B} / \hat{A}; \hat{\sigma}} =
\opn{Tr}_{\hat{B} / \hat{A}; \rho} =
\opn{Tr}_{\hat{B} / \hat{A}; \hat{\tau}} . \]
But $q_{\hat{A} / A}$ is injective, so the claim is proved. Our arguments
also imply properties (i), (ii) and (iii).

Let us now prove that $\opn{Tr}_{B/A}$ is continuous. First assume that
$\opn{res.dim} f \leq 1$. Then
$\cal{K}(B)$, being a cofinite type ST $B$-module, actually has the fine
$A$-module topology. (cf.\ \cite{Ye1} Def.\ 3.3 or Def.\ 3.2.1 (b.ii)). Since
$\opn{Tr}_{B/A} : \cal{K}(B) \ar \cal{K}(A)$ is $A$-linear, it is
continuous. Now assume $\opn{res.dim} f = n > 1$.
Consider the prime ideal
$\frak{p} := \opn{Ker}(A \ar \kappa_{n-1}(B))$.
We can assume that $\frak{p} \neq \mx$, by replacing (if necessary) $A$
with $A [\sqbr{ t }]$, and sending $t$ to a parameter of
$\cal{O}_{n}(B)$. Thus $A / \frak{p}$ is a DVR and
$C := \lim_{\leftarrow i} (A / \frak{p}^{i})_{\frak{p}}$
is a BCA. The morphism $A \ar B$ factors into morphisms $A \ar C \ar B$,
both of $\opn{res.dim} < n$. By
induction $\opn{Tr}_{B/C}$ and $\opn{Tr}_{C/A}$ are continuous, and
$\opn{Tr}_{B/A} = \opn{Tr}_{B/C}\ \circ \opn{Tr}_{C/A}$.

Finally to prove (iv), take a coefficient field $\sigma: K \ar A$. Then
\[ \Psi^{B}_{f \circ \sigma} : \cal{K}(B) \ar
\opn{Dual}_{f \circ \sigma} B =
\opn{Hom}_{K}^{\mathrm{cont}}(B, \omega(K)) \]
is bijective. On the other hand, one easily sees that
\[ \opn{Hom}_{A}^{\mathrm{cont}}(B, \cal{K}(A)) \ar
\opn{Hom}_{K}^{\mathrm{cont}}(B, \omega(K)) \]
is injective, so
$\cal{K}(B) \iso \opn{Hom}_{A}^{\mathrm{cont}}(B, \cal{K}(A))$.
\end{pf}

\begin{rem} \label{rem7.1}
Suppose $A, B$ are BCAs,
$f: A \ar B$ is a continuous $k$-algebra homomorphism, and $M$ is a
torsion type ST $A$-module. For instance,
$A, B$ could be any complete local $k$-algebras which are residually
finitely generated over $k$, $f$ could be any
local homomorphism, and $M$ any $0$-dimensional $A$-module.
If $M$ has finite length, define
\[ f_{\#} M := \opn{Dual}_{B} (B \otimes_{A} \opn{Dual}_{A} M). \]
Otherwise $M = \lim_{\alpha \ar} M_{\alpha}$ where each $M_{\alpha}$ has
finite length, and we set
$f_{\#} M := \lim_{\alpha \ar} f_{\#} M_{\alpha}$.
This gives a functor
$f_{\#} : \mathsf{STMod}_{\mathrm{tors}}(A) \ar
\mathsf{STMod}_{\mathrm{tors}}(B)$. Note that
$f_{\#} \cal{K}(A) = \cal{K}(B)$. If $f$ is a morphism in $\mathsf{BCA}(k)$,
the trace map $\opn{Tr}_{f}: \cal{K}(B) \ar \cal{K}(A)$ defines a trace map
$\opn{Tr}_{f} : f_{\#} M \ar M$ for any $M$.
The collection of data
$(\mathsf{STMod}_{\mathrm{tors}}(A), f_{\#})$ is a realization (and
generalization) of Lipman's pseudofunctor on $0$-dimensional modules;
cf.\ \cite{Hg}.
\end{rem}


\section{Duals of Continuous Differential Operators}

In this section we consider a continuous differential operator
$D : M \ar N$, and construct a dual operator
$\opn{Dual}_{A}(D) : \opn{Dual}_{A} N \ar \opn{Dual}_{A} M$.
The idea is to use the right $\cal{D}(K)$-module structure of
$\omega(K)$, for a TLF $K$.

Let $A$ be a local BCA with maximal ideal $\mx$, and let
$\sigma : K \ar A$ be a pseudo coefficient field. Given two finite type
ST $A$-modules $M,N$, choose $\mx$-filtered $K$-bases
$\ul{x} = (x_{0}, x_{1}, \ldots)$ and
$\ul{y} = (y_{0}, y_{1}, \ldots)$ for $M$ and $N$, respectively (cf.\ Def.\
\ref{def6.3}).
Suppose $D : M \ar N$ is a continuous DO over $A$ relative to $k$. For
$i,j \geq 0$ let
$D_{ij} : K \ar K$ be the functions such that, for $\lambda \in K$,
\[ D(\sigma(\lambda) x_{i}) = \sum_{j} \sigma(D_{ij}(\lambda)) y_{j} . \]
Then, just like in Lemma \ref{lem6.2}, $D_{ij} \in \cal{D}(K)$.

\begin{dfn} \label{dfn8.1}
Let
$\opn{Dual}_{\sigma}(D) : \opn{Dual}_{\sigma} N \ar \opn{Dual}_{\sigma} M$
be the function taking
$\phi \in \opn{Dual}_{\sigma} N$ to
\[ \opn{Dual}_{\sigma}(D) (\phi) :
\sum_{i} \sigma(\lambda_{i}) x_{i} \mapsto
\sum_{i,j} \lambda_{i} (\phi(y_{j}) * D_{ij})  . \]
\end{dfn}

There is no reference in the notation ``$\opn{Dual}_{\sigma}(D)$'' to the bases
$\ul{x}, \ul{y}$. This is not an oversight - as we shall see, this function
is independent of the bases. First, another definition:

\begin{dfn} \label{dfn8.0}
Let $M$ be a ST $A$-module (not necessarily of finite type).
Define the residue pairing to be
\begin{eqnarray*}
\langle -,- \rangle_{A/K}^{M} & : & M \times \opn{Dual}_{A} M \ar \omega(K) \\
\langle x, \phi \rangle_{A/K}^{M} & = & \opn{Res}_{A/K}(\phi(x))
\end{eqnarray*}
where $\opn{Res}_{A/K}$ is as in Def.\ \ref{dfn7.2}.
\end{dfn}

\begin{rem} \label{rem8.1}
Suppose $K$ is discrete (i.e.\ $\opn{dim} K = 0$) and $M$ is a finite type
or a cofinite type ST
$A$-module. Then the topology on $M$ is $K$-linear (cf.\ \cite{Ye1} Prop.\
3.2.5). As a topological vector space over $K$, $M$ is strongly
reflexive, in the sense of \cite{Ko} \S 13.3. One can show that the strong
$\opn{Hom}_{K}$ topology on
$\opn{Dual}_{A} M \cong \opn{Hom}_{K}^{\mathrm{cont}}(M, \omega(K))$ coincides
with the fine $A$-module topology on it. Hence  \linebreak
$\langle -,- \rangle_{A/K}^{M}$ is a perfect pairing also from the point
of view of \cite{Ko}.
\end{rem}

\begin{lem} \label{lem8.2} \mbox{ }

\begin{enumerate}
\rmitem{a} Suppose $\opn{ord}_{K}(D) = 0$, i.e.\ $D$ is $K$-linear. Then
$\opn{Dual}_{\sigma}(D) (\phi) = \phi \circ D$ for all
$\phi \in \opn{Dual}_{\sigma} N$.

\rmitem{b} Suppose $k \ar K$ is a morphism in $\mathsf{BCA}(k)$. Then
for all $\phi \in \opn{Dual}_{\sigma} N$,
\[ \opn{Res}_{K/k} \circ \opn{Dual}_{\sigma}(D) (\phi) =
\opn{Res}_{K/k} \circ \phi \circ D . \]
In other words, $\opn{Dual}_{\sigma}(D)$ is adjoint to $D$ with respect to the
the residue pairings
$\langle -,- \rangle_{A/k}^{M}$ and
$\langle -,- \rangle_{A/k}^{N}$.
\end{enumerate}
\end{lem}

\begin{pf}
One has $D_{ij} = \mu_{ij} \in K \subset \cal{D}(K)$, where
$D(x_{i}) = \sum_{j} \sigma(\mu_{ij}) y_{j}$. Now simply plug this into the
definition of $\opn{Dual}_{\sigma}(D)$.

\medskip \noindent (b)\ Say
$\phi(y_{j}) = \alpha_{j} \in \omega(K)$. Given
$x = \sum_{i} \sigma(\lambda_{i}) x_{i} \in M$, with $\lambda_{i} \in K$,
we have by the definition of the DOs $D_{ij}$:
\[ D(x) = \sum_{i,j} \sigma(D_{ij} * \lambda_{i}) y_{j}\ , \]
so
\begin{eqnarray*}
\lefteqn{ \langle D(x), \phi \rangle_{A/k}^{N} =
\opn{Res}_{K/k} \circ \phi \circ D (x) }\\
& & = \opn{Res}_{K/k} (\sum_{i,j}  (D_{ij} * \lambda_{i}) \alpha_{j}) =
\sum_{i,j} \langle D_{ij} * \lambda_{i}, \alpha_{j} \rangle_{K/k}.
\end{eqnarray*}
On the other hand, by the definition of $\opn{Dual}_{\sigma}(D)$,
\begin{eqnarray*}
\lefteqn{ \langle x, \opn{Dual}_{\sigma}(D)(\phi) \rangle_{A/k}^{M} =
\opn{Res}_{K/k} \circ (\opn{Dual}_{\sigma}(D) (\phi)) (x) } \\
& & = \opn{Res}_{K/k} (\sum_{i,j} \lambda_{i} (\alpha_{j} * D_{ij})) =
\sum_{i,j} \langle \lambda_{i}, \alpha_{j} * D_{ij} \rangle_{K/k}.
\end{eqnarray*}
Now use Thm.\ \ref{thm5.1}.
\end{pf}

\begin{lem} \label{lem8.3} \mbox{ }
\begin{enumerate}
\rmitem{a} $\opn{Dual}_{\sigma}(D) : \opn{Dual}_{\sigma} N \ar
\opn{Dual}_{\sigma} M$
is a continuous DO over $A$, relative to $k$, of order $\leq \opn{ord}_{A}(D)$.
It is independent of the $\mx$-filtered $K$-bases $\ul{x}, \ul{y}$.

\rmitem{b} Let $\mx \subset A$ be the maximal ideal. Suppose
$\sigma' : K \ar A$ is another pseudo coefficient field, s.t.\
$\sigma' \equiv \sigma\ (\opn{mod} \mx)$. Then
\[ \opn{Dual}_{\sigma'}(D) =  \Psi_{\sigma, \sigma'}^{M} \circ
\opn{Dual}_{\sigma}(D) \circ \Psi_{\sigma', \sigma}^{N} . \]

\rmitem{c} Suppose $\tau : L \ar A$ is another pseudo coefficient field, and
$f : K \ar L$ is a (finite) morphism in $\mathsf{BCA}(k)$, s.t.\
$\sigma = \tau \circ f$. Then for each $\phi \in \opn{Dual}_{\tau} N$,
\[ \opn{Dual}_{\sigma}(D) (\opn{Tr}_{f} \circ \phi) = \opn{Tr}_{f} \circ
\opn{Dual}_{\tau}(D) (\phi)  . \]
\end{enumerate}
\end{lem}

\begin{pf}
The proof resembles that of Prop.\  \ref{prop6.2}. Choose an intensification
homomorphism $u : K \ar \hat{K}$
such that $k \ar \hat{K}$ is a morphism. Let
$\hat{A} := A \otimes_{K}^{(\wedge)} \hat{K}$ and $v : A \ar \hat{A}$.
Replacing $A$ with each of the localizations $\hat{A}_{\hat{\mx}}$,
$\hat{\mx} \in \opn{Max} \hat{A}$, allows us to assume that $k \ar K$ is
itself a morphism in $\mathsf{BCA}(k)$. By Lemma \ref{lem8.2} (b) we see
that $\opn{Dual}_{\sigma}(D)$ is the adjoint of  $D$ w.r.t.\
the residue pairings $\langle -,- \rangle_{A/k}^{M}$ and
$\langle -,- \rangle_{A/k}^{N}$, so in particular it is independent of the
$\mx$-filtered $K$-bases $\ul{x}, \ul{y}$. It also follows
that for any $a \in A$,
\[ [\opn{Dual}_{\sigma}(D), a] = - \opn{Dual}_{\sigma}([D,a]) :
\opn{Dual}_{\sigma} M \ar \opn{Dual}_{\sigma} N, \]
bounding the order of the operator $\opn{Dual}_{\sigma}(D)$.
Here ``$[-,-]$'' denotes the commutator.
Parts (b),(c) of the present lemma are similarly proved, using Lemma
\ref{lem6.3}.

As for the continuity of $\opn{Dual}_{\sigma}(D)$, it can be deduced from the
fact that it is a linear combination of the continuous operators $D_{ij}$
appearing in its definition.
\end{pf}

The ST $A$-module $\cal{K}(A)$ is separated. Therefore for any ST $A$-module
$M$, the canonical surjection $M \surj M^{\opn{sep}}$ induces an isomorphism
$\opn{Dual}_{A} M^{\opn{sep}} \iso \opn{Dual}_{A} M$.
Here is the third main result of the paper:

\begin{thm} \label{thm8.1} \rom{(Duals of Continuous DOs)}\
Let $A$ be a  BCA over $k$. Let $M$ and $N$ be ST $A$-modules with the fine
topologies, and let $D : M \ar N$ be a continuous DO over $A$ relative to
$k$. Then there is a unique function
\[ \opn{Dual}_{A}(D) : \opn{Dual}_{A} N \ar \opn{Dual}_{A} M, \]
satisfying the conditions below:

\begin{enumerate}
\rmitem{i} $\opn{Dual}_{A}(D) : \opn{Dual}_{A} N \ar \opn{Dual}_{A} M$ is a
continuous DO over $A$ relative to $k$, of order $\leq \opn{ord}_{A}(D)$.

\rmitem{ii} (Transitivity)\ if $E : N \ar P$ is another such operator,
then $\opn{Dual}_{A}(E \circ D) = \opn{Dual}_{A}(D) \circ \opn{Dual}_{A}(E)$.

\rmitem{iii} (Linearity)\
if $D$ is $A$-linear, then $\opn{Dual}_{A}(D)$ is the homomorphism
$\phi \mapsto \phi \circ D$, for
$\phi \in \opn{Dual}_{A} N = \opn{Hom}^{\mathrm{cont}}_{A}(N, \cal{K}(A))$.

\rmitem{iv} (Base change)\ let $v : A \ar \hat{A}$ be an intensification
homomorphism, and let
$\hat{D} : (\hat{A} \otimes_{A} M)^{\opn{sep}}
\ar (\hat{A} \otimes_{A} N)^{\opn{sep}}$ be the unique extension of $D$. Then
\[ \opn{Dual}_{\hat{A}}(\hat{D}) \circ q^{N}_{v} =
q^{M}_{v} \circ \opn{Dual}_{A}(D), \]
where $q^{M}_{v}, q^{N}_{v}$ are the homomorphisms of Prop.\ \ref{prop7.1}.

\rmitem{v} Assume $\sigma: K \ar A$ is a morphism in $\mathsf{BCA}(k)$ s.t.\
$D$ is $K$-linear. Then $\opn{Dual}_{A}(D)$ is the adjoint to $D$ w.r.t.\ the
residue pairings
$\langle - , - \rangle^{M}_{A/K}$ and
$\langle - , - \rangle^{N}_{A/K}$.

\rmitem{vi} Suppose $A$ is local and $M,N$ are finite type ST
$A$-modules.
Given a pseudo coefficient field $\sigma : K \ar A$, one has
\[ \Psi_{\sigma}^{M} \circ \opn{Dual}_{A}(D) = \opn{Dual}_{\sigma} (D)
\circ \Psi_{\sigma}^{N} . \]
Here $\Psi_{\sigma}^{M}, \Psi_{\sigma}^{N}$ are the isomorphisms of Thm.\
\ref{thm6.1}, and $\opn{Dual}_{\sigma} (D)$ is the function defined in Def.\
\ref{dfn8.1}.
\end{enumerate}
\end{thm}

\begin{rem}
Trivially, the category $\mathsf{Mod}(A)$ of
$A$-modules and $A$-linear homomorphisms, and the category
$\mathsf{STMod}_{\opn{fine}}(A)$ of ST $A$-modules with fine topologies
and continuous
$A$-linear homomorphisms, are equivalent (under the functor
$\opn{untop} : \mathsf{STMod}(A) \ar \mathsf{Mod}(A)$ which forgets
the topology). However, if we take the same classes of objects, but
enlarge the set of morphisms between two objects to be DOs and continuous
DOs, respectively, these new categories are no longer equivalent. This is
so at least when $\opn{char} k = 0$ and $\opn{res.dim} A \geq 1$.
Our results are valid only for continuous DOs.
\end{rem}

\begin{pf}
Let $M,N$  be finite type ST $A$-modules. Using Lemma \ref{lem8.3}, and
proceeding just like in the proofs of Theorems \ref{thm6.1} and
\ref{thm7.2}, we arrive at a function $\opn{Dual}_{A}(D)$ which satisfies
conditions (i)-(iv), (vi).
As for condition (v), after a base change $K \ar \hat{K}$ we reduce to the
case when $k \ar A$ is a morphism. Now we can use Lemma \ref{lem8.2} (b).

Now let $M,N$ be ST $A$-modules with fine topologies. After possibly
applying $(-)^{\opn{sep}}$ to these modules, we may assume they are separated.
Choose an isomorphism $M \cong \lim_{\alpha \ar} M_{\alpha}$, with the
$M_{\alpha}$ modules of finite type. Let $N_{\alpha}$ be the $A$-module
$A \cdot D(M_{\alpha}) \subset N$, endowed with the fine topology.
Say $d = \opn{ord}_{A}(D)$. Because
$\cal{P}_{A/k}^{d,\opn{sep}}(M_{\alpha})$ is a finite type ST $A$-module
(by Prop.\ \ref{prop4.3}), $N_{\alpha}$ is of finite type, and
$D_{\alpha} := D|_{M_{\alpha}} : M_{\alpha} \ar N_{\alpha}$ is continuous.
Let $\psi : \lim_{\alpha \ar} N_{\alpha} \ar N$ be the inclusion, and set
\[ \opn{Dual}_{A}(D) :=  (\lim_{\leftarrow \alpha}
\opn{Dual}_{A}(D_{\alpha}) ) \circ \opn{Dual}_{A}(\psi) . \]
Since the functor $\opn{Dual}_{A}$ sends $\lim_{\ar}$ to $\lim_{\leftarrow}$
(cf.\ Lemma \ref{lem1.1} (4)),
this extended
definiton of $\opn{Dual}_{A}(D)$ satisfies all the conditions of the theorem.
\end{pf}

Occasionally we shall abbreviate $\opn{Dual}_{A} M$ to $M^{\vee}$, and
$\opn{Dual}_{A}(D)$ to $D^{\vee}$.

\begin{cor} \label{cor8.2} \mbox{ }

\begin{enumerate}
\rmitem{a} Let $M,N$ be each either finite type or cofinite type ST
$A$-modules, and let
$D \in \opn{Diff}_{A/k}^{\opn{cont}}(M,N)$. Then under the canonical
isomorphisms $M \iso M^{\vee \vee}$ and $N \iso N^{\vee \vee}$, one has
$D \mapsto D^{\vee \vee}$.

\rmitem{b} With $M,N$ as in \rom{(a)}, the map
$\opn{Diff}_{A/k}^{\opn{cont}}(M,N) \ar
\opn{Diff}_{A/k}^{\opn{cont}}(M^{\vee},N^{\vee})$,
$D \mapsto D^{\vee}$, is an anti-isomorphism of filtered $A$-$A$-bimodules.
In particular,
$\cal{D}(A;$ \linebreak
$\cal{K}(A)) \cong \cal{D}(A)^{\circ}$
as filtered $k$-algebras.
\end{enumerate}
\end{cor}

\begin{pf}
(a)\ Using base change we can assume that $k \ar A$ is a morphism in
$\mathsf{BCA}(k)$. Then both $D$ and $D^{\vee \vee}$ are adjoints to
$D^{\vee}$ w.r.t.\ the residue pairings
$\langle -,- \rangle_{A/k}^{M}$ and
 $\langle -,- \rangle_{A/k}^{N}$.

\medskip \noindent (b)\
Immediate from part (a).
\end{pf}

Here are a couple of examples to illustrate the scope of our results:

\begin{exa} \label{exa8.2}
Suppose $A$ is a noetherian, local, residually finitely generated
$k$-algebra. Let $I$ be an injective hull of the residue field
$A / \frak{m}$. Then $I$ is (non-canonically) a right $\cal{D}(A)$-module,
and moreover
$\opn{Diff}_{A/k}(I,I) \cong \cal{D}(\hat{A})^{\circ}$, where
$\hat{A}$ is the $\frak{m}$-adic completion. This is because $\hat{A}$ is
a BCA, there exists
an isomorphism of $\hat{A}$-modules $I \cong \cal{K}(\hat{A})$, and any
DO $I \ar I$ is automatically continuous for the $\frak{m}$-adic topology.
\end{exa}

\begin{exa} \label{exa8.3}
Let $A$ be a BCA.
Suppose $M^{\bdot}$ is a bounded complex with each $M^{q}$ a finite type
ST $A$-module, and $D: M^{q} \ar M^{q+1}$ a continuous DO (for instance,
$M^{\bdot} = \Omega^{\bdot, \mathrm{sep}}_{A/k}$). Then
$\opn{Dual}_{A} M^{\bdot}$ is also a complex (of cofinite type modules), and
a standard spectral sequence argument shows that the homomorphism of
complexes
\[ M^{\bdot} \ar \opn{Dual}_{A} \opn{Dual}_{A} M^{\bdot} \]
(in the abelian category of untopologized $k$-modules)
is a quasi-isomorphism.
\end{exa}

\begin{question} \label{que8.1}
In the example above, suppose the complex $M^{\bdot}$ is acyclic. Is the
same true of the dual complex $\opn{Dual}_{A} M^{\bdot}$? A slight variation
is: suppose
$\opn{rank}_{k} \mathrm{H}^{q} M^{\bdot}$ \linebreak
$< \infty$ for all $q$. Is the same
true for $\opn{Dual}_{A} M^{\bdot}$?
\end{question}

\begin{cor} \label{cor8.3}
Let $f: A \ar B$ be a morphism in $\mathsf{BCA}(k)$, let $M$ (resp.\ $N$)
be a ST $A$-module (resp.\ $B$-module) with the fine topology, and let
$D \in \opn{Diff}^{\mathrm{cont}}_{A/k}(M,N)$. Then there is a DO
\[ \opn{Dual}_{B/A}(D) = \opn{Dual}_{f}(D) : \opn{Dual}_{B} N \ar
\opn{Dual}_{A} M  . \]
The asignment $D \mapsto \opn{Dual}_{f}(D)$ satisfies the obvious
generalizations of conditions (i)-(v) of Thm. \ref{thm8.1}. For instance
(iii): if $D$ is $A$-linear, then
$\opn{Dual}_{f}(D)(\phi) = \opn{Tr}_{B/A} \circ \phi \circ D$.
\end{cor}

\begin{pf}
We may assume that $M$ is a finite type ST $A$-module, and that $N$ is
separated. So $D$ factors into
$M \exar{ \mathrm{d}^{n}_{M} } \cal{P}^{n,\mathrm{sep}}_{A/k}(M)
\exar{ \phi } N$, with
$\phi \in \opn{Hom}_{A}^{\mathrm{cont}}(\cal{P}^{n,\mathrm{sep}}_{A/k}(M), N)$
and $n \geq \opn{ord}_{A}(D)$. Let
$\phi^{\vee} : \opn{Dual}_{B} N \ar \opn{Dual}_{A}
\cal{P}^{n,\mathrm{sep}}_{A/k}(M)$
be the homomorphism
$\psi \mapsto \opn{Tr}_{B/A} \circ \psi \circ \phi$, for
$\psi \in \opn{Dual}_{B} N =\opn{Hom}_{B}^{\mathrm{cont}}(N, \cal{K}(B))$.
Define
$\opn{Dual}_{f}(D) := \opn{Dual}_{A}(\mathrm{d}^{n}_{M}) \circ \phi^{\vee}$.
The transitivity and uniqueness properties follow from base change
and the uniqueness of adjoints.
\end{pf}

\begin{exa} \label{exa8.1}
If $f : A \ar B$ is a morphism of BCAs,
the trace map
$\opn{Tr}_{B/A}: \cal{K}(B) \ar \cal{K}(A)$ and
the continuous DGA homomorphism
$\Omega^{\bdot, \mathrm{sep}}_{A/k} \ar \Omega^{\bdot, \mathrm{sep}}_{B/k}$
induce a map
$\opn{Tr}_{B/A}: \opn{Dual}_{B} \Omega^{\bdot, \mathrm{sep}}_{B/k} \ar
\opn{Dual}_{A} \Omega^{\bdot, \mathrm{sep}}_{A/k}$,
which by the corollary is a homomorphism of complexes. This fact is important
for the construction of the De Rham - residue double complex in \cite{Ye2}.
\end{exa}


\end{document}